\documentclass{article}

\usepackage{amssymb}
\usepackage{color}
\usepackage{arxiv}
\usepackage{dsfont}
\usepackage{amsmath}
\usepackage{amsthm}
\usepackage{amsfonts}
\newtheorem{theorem}{Theorem}
\newtheorem{corollary}{Corollary}

\usepackage{float}
\usepackage[utf8]{inputenc} 
\usepackage[T1]{fontenc}    
\usepackage{hyperref}       
\usepackage{url}            
\usepackage{booktabs}       
\usepackage{amsfonts}       
\usepackage{nicefrac}       
\usepackage{microtype}      
\usepackage{lipsum}		
\usepackage{graphicx}
\usepackage{natbib}
\usepackage{doi}
\usepackage{multirow}
\usepackage{rotating}
\usepackage{float}
\hypersetup{                    
    colorlinks=true,                
    breaklinks=true,                
    urlcolor= blue,                 
    linkcolor= blue,                
    citecolor= green                
    }

\title{Bayesian composite confidence interval for the tail index under randomly right-censored data}

\author{ \href{https://orcid.org/0000-0003-3855-9486}{\includegraphics[scale=0.06]{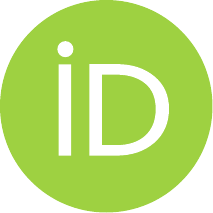}\hspace{1mm}Abdelkader Ameraoui}\thanks{Corresponding author: A. Ameraoui. \quad Address: National Higher School of Mathematics. Scientific and Technology Hub of Sidi Abdellah, Po. Box 75, Mahelma 16093, Algiers (Algeria) -- E-mail: \texttt{aameraoui@nhsm.edu.dz}.} \\
	National Higher School of Mathematics\\ Scientific and Technology Hub of Sidi Abdellah\\
 Po. Box 75, Mahelma 16093, Algiers (Algeria) \\
	\texttt{aameraoui@nhsm.edu.dz} \\
	\And
	\href{https://dupuy.perso.math.cnrs.fr/}{\includegraphics[scale=0.04]{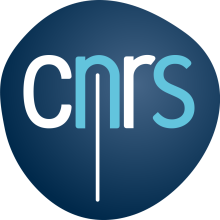}\hspace{1mm}Jean-Fran\c cois ~Dupuy} \\
	Univ. de Rennes, INSA de Rennes, CNRS,\\ IRMAR - UMR 6625, F-35000 Rennes, (France) \\
	\texttt{Jean-Francois.Dupuy@insa-rennes.fr } \\
 \And 
 \href{}{\hspace{1mm}Kamal ~Boukhetala} \\
	Faculty of Mathematics, University of Science and Technology Houari Boumediene,\\
Po. Box 32 El-Alia 16111, Algeria\\
	\texttt{kboukhetala@usthb.dz } \\
	 \AND
}

\renewcommand{\shorttitle}{Bayesian confidence interval for tail index}

\hypersetup{
pdftitle={Bayesian composite confidence interval for tail index under right censorship},
pdfsubject={math.ST},
pdfauthor={Abdelkader~Ameraoui, Jean-Fran\c cois~Dupuy},
pdfkeywords={Bayesian inference, Composite likelihood, tail index, },
}
\date{}
\begin{document}
\maketitle
\begin{abstract}
Bayesian composite likelihood estimation of the tail index of a heavy-tailed distribution is addressed when
data are randomly right-censored. Maximum a posteriori and mean posterior estimators
are constructed under Jeffrey's prior distribution of the tail index. Based on asymptotic results, some confidence regions (CR) for the tail index are constructed using posterior distribution and log-posterior ratio statistic. The proposed confidence regions are investigated via Finite-sample simulations. Finally, the proposed confidence regions are outperformed through two real datasets.
\end{abstract}


\section{Introduction}

In modern statistical literature, the occurrence of extreme events generates catastrophic scenarios. Extreme value theory (EVT) is a branch of statistics that focuses on modeling the risk associated with rare and extreme events. These events are often outliers, significantly deviating from the central tendency of a distribution. 
In this context, the tail index estimation is a crucial issue in extreme value theory, that measures the thickness of the tail of a probability distribution function and plays a major role for evaluating the risk of occurrence of  extremes events. Extreme events arise in a wide variety of domains and a vast literature has been dedicated to this topic. Recent overviews can be found in the monographs \cite{Bei:04}, \cite{Emb:01} and \cite{dHan01}.

A relatively recent statistical approaches, alternative to the classical maximum likelihood method, are a non-parametric analogue, called empirical likelihood \cite[see, e.g.,][]{owen01,owen02,owen03},  or more recently \cite{owen04} and the composite likelihood (or the weighted likelihood) estimation approaches that can be assimilated to a nonparametric version of the Wilks theorem \cite[see, e.g. ][]{Cas:24}. The likelihood ratio tests statistics based on the empirical likelihood entail a parametric properties of the 
likelihoods such an asymptotic $\chi^2$ distribution allowing the construction of more precise confidence regions. Indeed, this approach remains valid, and can be applied in a regression models \cite[see, e.g.,][]{chen01}, censoring framework estimation \cite[see, e.g.,][]{he01,liang01}, but also in a Bayesian context \cite[see, e.g.,][]{laz01}.

Composite and empirical likelihood estimation provide an alternative to asymptotic confidence interval based on the asymptotic normality  approach. In the context of extreme values
analysis, weighted likelihood estimators have been investigated in many works as \cite{Stn:23}, \cite{li01}, \cite{qi01}, or \cite{peng01}. In risk measurement context, we cite \cite{bay01}, \cite{zhen01}  \cite{peng02}, but to the best of our knowledge, no confidence interval for the tail index has been proposed when the data are randomly right censored and covers of heavy-tailed distribution. The present work intends to construct confidence regions of the tail index based on Bayesian posterior distribution and composite likelihood -both in the case of complete and randomly right censored data - to fill in this gap.

The paper is organized as follows. In Section \ref{sec2} and \ref{sec3}, we construct our confidence regions based on Bayesian composite likelihood approach and consider asymptotic results. Proofs are deferred to an appendix. Section \ref{sec4} reports the results of a comprehensive simulation study. Finite-sample performance of the proposed confidence regions is assessed via simulations. Finally, we illustrate our methodology on a real dataset of global costly natural disasters (for complete data) and Australian AIDS survival data (when randomly right censorship is considered) in Section \ref{sec5}. A discussion and some perspectives are given in Section \ref{sec6}.

\section{Confidence interval for the tail index} \label{sec2}
Let $F$ be the cumulative distribution function (cdf) of some non-negative random variable $X$. We assume that $F$ is heavy-tailed, that is, there exists a constant $\alpha>0$ such that
\begin{equation}\label{eq.1}
1-F(x)=x^{-\alpha}\ell(x),
\end{equation} 
where $\ell$ is a slowly varying function at infinity:
\[
\lim_{x \rightarrow \infty} \frac{\ell(tx)}{\ell(x)}=1 \qquad\mbox{for all } t>0.   
\]
 
If (\ref{eq.1}) holds, we have:
\[
\lim_{x\rightarrow \infty} \frac{1-F(tx)}{1-F(x)}=t^{-\alpha} \qquad\mbox{for all } t>0
\] 
and we say that $\bar{F}=1-F$ is regularly varying at infinity with tail index $\alpha$, which we denote by $\bar{F} \in \mathcal{R}_{-\alpha}$. The positive number $\gamma:=\alpha^{-1}$ is called the extreme value index (EVI) of $F$. The conditions above amount to assuming that the  distribution function $F$  is in the max-domain of attraction of a Fr\'echet distribution. Such distribution functions are useful in practice for investigating phenomena where exceptional values have a significant occurrence frequency. 

Several estimators have been proposed for the tail index $\alpha$, or equivalently, for the EVI $\gamma$ \citep[see, e.g.,][]{Ber:01,Ber:02}. Let $X_1,X_2,\ldots ,X_n$ be independent and identically distributed (iid thereafter) random variables with common cdf $F$. Let $k\in \{2,\ldots,n\}$ and $X_{n,1}\leq X_{n,2} \leq \ldots \leq X_{n,n}$ be the order statistics of the sample $X_1,X_2,\ldots ,X_n$. Hill estimator is defined as
\begin{equation}\label{eq.2}
H(k):=\frac{1}{k}\sum_{i=1}^{k}{\log(X_{n,n-i+1})-\log(X_{n,n-k})}.
\end{equation}
Consistency of this estimator was proved in \cite{hill01} under the regular variation condition (\ref{eq.1}). Its asymptotic normality was further established under an additional condition known as the  second-order regular variation condition  \citep[see,][]{dHan01}, where we assume that there exists a function $A(t) \to 0$ as $t \to \infty$, such that :
\begin{equation}
    \lim_{t \to  \infty} \dfrac{\frac{U(tx)}{U(t)}-x^{-\alpha}}{A(t)}=x^{\alpha} \frac{x^\beta-1}{\beta}
    \label{eq1.5}
\end{equation}
for all $x > 0$, $\beta \leq 0$, and $U(t):=(1/(1-F))^{\leftarrow}(t)=\inf\{x:F(x)\leq 1-\frac{1}{t}$\}  (We say that the function $U$ is of regular variation of second order, and $A$ is regularly varying function with index $\beta$).

In this section, we introduce our methodologies by assuming complete data set with commun distribution $F$ has the simplified form of (\ref{eq.1}), given by :
\begin{equation}
\label{eq2.0}
    1-F(x)=c x^{-\alpha}, \mbox{ for } x>u, c>0
\end{equation}
where $u$ is a suitable threshold, and let $\delta_i= \mathds{1}_{\left\lbrace X_i > u \right\rbrace }$. Then, the likelihood function for the framework $\{\delta_i, \max(X_i,u)\}_{i=1}^n$ is 
\[
L(\alpha , c) = \prod_{i=1}^n \left( c \alpha X_i^{-1-\alpha}\right)^{\delta_i} \left( 1- c  u^{-\alpha}\right)^{1-\delta_i}
\]
In this paper, we let $u=X_{n,n-k}$, where $k=k(n)$ is a sequence satisfying
\begin{equation}
    \label{eq2.1}
    k \to \infty \mbox{ and } \frac{k}{n} \to 0 \mbox{ when } n \to \infty.
\end{equation}
Then, the likelihood function above becomes
\begin{equation}
    \label{eq2.2}
 L(\alpha , c) = \prod_{i=1}^n \left( c \alpha X_i^{-1-\alpha}\right)^{\delta_i} \left( 1- c  X_{n,n-k}^{-\alpha}\right)^{1-\delta_i}   
\end{equation}
Now, we focus on defining more general form of (\ref{eq2.2}), called the composite-likelihood (CL) functions, that have emerged gradually through  \cite{Lind:88}. This model offers several advantages over traditional likelihood methods, particularly when dealing with complex data or models are defined under conditional assumptions. Our methodology is inspired by the Bayes confidence interval construction in \cite{Efr:93}, and motivated by the results of \cite{laz01}, and \cite{Pau:12}, \cite{Rib:12} where the application of these approaches on spatial extreme data are discussed and showed that the posterior credible intervals yield appropriate empirical coverage rates. More recently \cite{chan:17} discussed the performance of the Bayesian composite likelihood estimation when max-stable processes is considered. \cite{Stn:23} proposed parametric estimation of extreme quantiles of a distribution using weighted composite log-likelihood approach and showed that including weights in the composite log-likelihood function can reduce the sensitivity of estimates to small changes in the threshold. Confidence intervals are also constructed by inverting a test statistic calibrated via parametric bootstrapping. 

Based in (\ref{eq2.0}), we define the composite-likelihood function

\begin{equation}
    \label{eq2.22}
 L_{CL}(\alpha , c) = \prod_{i=1}^n \left( c \alpha X_i^{-1-\alpha}\right)^{nw_i\delta_i} \left( 1- c  X_{n,n-k}^{-\alpha}\right)^{nw_i(1-\delta_i)}   
\end{equation}
where $w=(w_1,w_2, \ldots,w_n)$ is a vector of non-negative weights with summation equal to one.

One can consider from (\ref{eq2.0}), the Jeffrey prior for $\alpha$ and $c$ using the Fisher information matrix given by
$$
J(\alpha,c)=\left[
\begin{array}{ll}
    \frac{cu^{-\alpha}}{\alpha^2} + \frac{cu^{-\alpha}\log ^2 u}{1-cu^{-\alpha}}& -\frac{u^{-\alpha}\log u}{1-cu^{-\alpha}} \\
     -\frac{u^{-\alpha}\log u}{1-cu^{-\alpha}} & -\frac{u^{-\alpha}}{c(1-cu^{-\alpha})}
\end{array}
\right]
$$
The posterior distribution of ($\alpha,c$) can be proportionally obtained with the composite likelihood in (\ref{eq2.22}) substituting for the full likelihood, as

\begin{align}
  \pi(\alpha, c\vert X_1,\ldots,X_n) &\propto  \pi(\alpha,c) \times L_{CL}(\alpha , c) \nonumber \\
  &\propto \frac{X_{n,n-k}^{-\alpha}}{\alpha} \left(1-cX_{n,n-k}^{-\alpha}\right)^{-1/2}\prod_{i=1}^n \left( c \alpha X_i^{-1-\alpha}\right)^{nw_i\delta_i} \left( 1- c  X_{n,n-k}^{-\alpha}\right)^{nw_i(1-\delta_i)}
  \label{eq2.3}
\end{align}

Under the posterior distribution in (\ref{eq2.3}), our approach consists on the  construction of Bayesian confidence regions of $\alpha$ using the data tilting method (\citep[see, e.g.][]{peng01}):

 First, for a given weights $w=(w_1,w_2, \ldots , w_n)$, such that $w_i \geq 0 \mbox{ for any } i=\overline{1,n}$ and $\sum_{i=1}^n w_i  =1$, we maximize the weighted log-posterior distribution given by :

\begin{align*}
    \ell_1(\alpha , c) &= \log \pi(\alpha, c\vert X_1,\ldots,X_n) \\
    &= -\alpha \log X_{n,n-k} - \log \alpha -\frac{1}{2} \log (1-c X_{n,n-k}^{-\alpha} ) + n\sum_{i=1}^n w_i \delta_i \log \left( c \alpha X_i^{-1-\alpha}\right)^{\delta_i} + n\sum_{i=1}^n w_i (1-\delta_i)\log \left( 1- c  X_{n,n-k}^{-\alpha}\right) 
\end{align*}

and we define $(\hat{\alpha}(w),\hat{c}(w))$ as the estimators of $(\alpha,c)$ solution of the equation $$(\hat{\alpha}(w),\hat{c}(w))=\arg \max_{\alpha>0,c>0} \ell_1 (\alpha , c)$$
This results in
$$
\begin{cases}
    \hat{\alpha}(w) = \dfrac{\sum_{i=1}^n w_i \delta_i -\frac{1}{n}}{\sum_{i=1}^n w_i \delta_i \log \frac{X_i}{X_{n,n-k}} + \frac{\log X_{n,n-k}}{n}} \\
    \hat{c}(w) = X_{n,n-k}^{\hat{\alpha}(w)} \frac{2n}{2n-1} \sum_{i=1}^n w_i \delta_i.
\end{cases}
$$
Let define the function $\mathbb{D}(w)=\sum_{i=1}^n w_i \log(n w_i)$, which is the measure of distance between $w$ and uniform distribution $w_i=1/n$. Our method consists in choosing the weights $w_i$ to minimize the distance $(2n)^{-1}\mathbb{L}(\alpha)=\mathbb{D}(w)$ according to the constraints
\begin{equation}
 w_i \geq 0, \qquad \sum_{i=1}^n w_i  =1, \qquad \sum_{i=1}^n w_i \delta_i \lbrace \log\left(\frac{X_{i}}{X_{n,n-k}} \right)-\frac{1}{\alpha} \rbrace =0.
 \label{eq2.10}
\end{equation}

Note that the constraint $\sum_{i=1}^n w_i \delta_i \lbrace \log\left(\frac{X_{i}}{X_{n,n-k}} \right)-\frac{1}{\alpha} \rbrace =0$ results from the fact that for all $i=1,\ldots, k$, $ y_i=~i~ \log\left(\frac{X_{n,n-i+1}}{X_{n,n-i}} \right) \overset{d}{=} \log\left(\frac{X_{n,n-i+1}}{X_{n,n-k}}\right)$,  are approximately $i.i.d$ exponential random variables with mean $\frac{1}{\alpha}$ (see e.g., \cite{weiss01}). Then, using standard Lagrange multipliers, we get  $w_i=w_i(\lambda_1, \lambda_2)$
\begin{equation}
= \begin{cases}
    \frac{1}{n} \exp (-1-\lambda_1) \quad \mbox{ if } \delta_i=0, \\
    \frac{1}{n} \exp \biggl\{ -1-\lambda_1  
 - \lambda_2 \biggl[ \log\left(\frac{X_{i}}{X_{n,n-k}} \right)-\frac{1}{\alpha} \biggr] \biggr\} \quad \mbox{ if } \delta_i=1,
    \end{cases}
    \label{eq2.110}
\end{equation}
where $\lambda_1, \lambda_2$ satisfy

\begin{equation}
 w_i(\lambda_1, \lambda_2) \geq 0, \qquad \sum_{i=1}^n w_i(\lambda_1, \lambda_2)  =1, \qquad \sum_{i=1}^n w_i(\lambda_1, \lambda_2) \delta_i \lbrace \log\left(\frac{X_{i}}{X_{n,n-k}} \right)-\frac{1}{\alpha} \rbrace =0.
 \label{eq2.100}
\end{equation}

\begin{theorem}\label{theo01}
Assume that the conditions in (\ref{eq1.5}) holds, when the sequence $k=k_n$, be such that $k \to \infty$ and $k/n \to 0$, $\sqrt{k} A(n/k) \to 0$, as $n \to \infty$.. Then it there exists a solution $(\lambda_1,\lambda_2)$ of (\ref{eq2.100}), such that
$$
\mathbb{L}(\alpha^0) \xrightarrow{\mathcal{D}} \chi^2_1,
$$
where $(\lambda_1,\lambda_2)=(\lambda_1(\alpha^0),\lambda_2(\alpha^0))$ in the definition of $\mathbb{L}(\alpha^0)$.

\end{theorem}
Therefore, and based on the above limit, an approximate 100($1-\theta$)\%  confidence region for $\alpha$ is
\begin{equation}
    I(\theta) =\{ \alpha : \mathbb{L}(\alpha) \leq q_\theta\}
    \label{eq2.12}
\end{equation}

where $q_\theta$ is the $\theta$-level critical point of $\chi^2_1$.

\section{Bayesian confidence region for the tail index under randomly right censored data}\label{sec3}
In this section, we address the construction of confidence region for the tail index $\alpha$ when data are randomly right-censored. Censoring commonly occurs in the analysis of event time data but also in extreme data analysis. For example, $X$ may represent the amount of claim related to an insurance policy, that can not be definitively evaluated when the data are collected, the information brought by this data is not completely available (or partially observed). An appropriate way to model this situation is to introduce a random variable $Y$ (called a censoring random variable) such that observations consist of pairs $(Z_i,\delta_i), 1\leq i \leq n$ where $Z_i=\min(X_i,Y_i)$, $\delta_i=\mathds{1}_{\left\lbrace X_i \leq Y_i \right\rbrace }$ and $\mathds{1}$ is the indicator function. Estimation of the EVI with censored data was considered in \cite{Wor:14}, \cite{Bei:07}, \cite{Bra:13}, \cite{Ein:01} and \cite{GN:11}. A general form of the EVI estimation in the context of censored data was proposed by \citeauthor{Ein:01} (\citeyear{Ein:01}) \cite{Ein:01} to estimate $\gamma$ by
\begin{equation}\label{eq.4}
\hat{\gamma}_Z^{(C)}(k)=\frac{\hat{\gamma}_Z(k)}{\hat{p}},
\end{equation}
where $\hat{\gamma}_Z(k)$ is any of the classical EVI estimators calculated on the censored observations $Z_1,\ldots,Z_n$ and $\hat{p}=\frac{1}{k}\sum_{i=1}^{k}\delta _{[n-i+1]}$ is the proportion  of uncensored values in the $k$ largest observations of $Z$. In this paper, we adopt a completely different approach and investigate Bayesian estimation of the tail index $\alpha:=\gamma^{-1}$. Bayesian estimation will allow us to incorporate \textit{a priori} knowledge about the data, and provides an alternative to frequentist methods.

In the context of extreme values analysis without censoring, Bayesian estimators have been investigated in \cite{Cab:11}, \cite{Coles:96}, \cite{Dieb:05}, \cite{Nasc:12}, \cite{Ber:02}. See also \cite{Bei:04} (chapter 11). A Bayesian estimator of the tail index $\alpha$ has been proposed by \cite{amer} when censoring is present, but to the best of our knowledge, no confidence interval of the tail index has been proposed in this context. 

In this section, we construct several Bayesian confidence regions for the tail index $\alpha$ in model \eqref{eq.1}. Bayesian estimation requires  specifying a prior distribution for the unknown parameter. Some asymptotic results on the the log posterior ratio are established. 

\textbf{Framework assumptions and  notations.} Let $X_1,X_2,\ldots,X_n$ (resp. $Y_1,Y_2,\ldots ,Y_n$) be $n$  iid copies of a non-negative  random variable $X$ with cdf $F$ (resp. $Y$ with cdf $G$). The probability density function of $X$ is denoted by $f$ (resp. of $Y$ is denoted $g$). We assume that $F$ (resp. $G$) is heavy-tailed with tail index $\alpha$ (resp. $\beta$), \textit{i.e.}, $\bar F \in \mathcal R_{-\alpha}$ (resp. $\bar G \in \mathcal R_{-\beta}$). We assume that $X$ and $Y$ are independent and that we observe the $n$ independent pairs
\begin{equation*}
(Z_i,\delta_i), 1\leq i \leq n,
\end{equation*}
where $Z_i=\min(X_i,Y_i)$ and $\delta_i=\mathds{1}_{\left\lbrace X_i \leq Y_i \right\rbrace }$. Let $Z_{n,1}\leq Z_{n,2} \leq \ldots \leq Z_{n,n}$ be the order statistics of  the sample  $(Z_1,Z_2,\ldots ,Z_n)$ and $\delta_{[n-i+1]}$ be the concomitant value of $\delta$ associated with $Z_{n,n-i+1}$. Let $H$ be the cdf of $Z$. Note that $H$ is also heavy-tailed and $\bar H \in \mathcal{R}_{-(\alpha +\beta )}$, by independence of $X$ and $Y$.

Let $E_{j,u}:=\frac{Z_{n,n-N_u+j}}{u}$, given $Z_{n,n-N_u+j}>u$, be the $j$-th relative excess over a threshold $u$ and $N_u$ denote the number of such excesses. For $j=1,\ldots,N_u$,  $E_{j,u}$ satisfies $\mathbb{P}(E_{j,u} >t ) \to t^{-\alpha}$ for $t>1$, as $u \to +\infty$. We can obtain the partial likelihood of $\alpha$ based on the sample $(E,\Delta)=(e_{j,u},\delta_{[n-N_u+j]})_{j=1,\ldots,N_u}$ (see e.g., \cite{Bei:07}):
\begin{eqnarray}\label{partlik}
\mathcal L_u(E,\Delta \vert \alpha)&=&\prod_{j=1}^{N_u}{\left[ \alpha {e_{j,u}}^{-(\alpha +1)}\right] ^{\delta_{[n-N_u+j]}}\left[ {e_{j,u}}^{-\alpha}\right] ^{1-\delta_{[n-N_u+j]}}} \nonumber \\
&=& \alpha ^{\sum_{j=1}^{N_u}{\delta_{[n-N_u+j]}}}\left( \prod_{j=1}^{N_u}{e_{j,u}}\right)^{-\alpha}\left( \prod_{j=1}^{N_u}{{e_{j,u}}^{-\delta_{[n-N_u+j]}}}\right).
\end{eqnarray}

Maximum likelihood estimator of $\alpha$ is given by
$$
\hat{\alpha}=\dfrac{\sum_{j=1}^{N_u} \delta_{[n-N_u+j]} }{\sum_{j=1}^{N_u} \log e_{j,u}}
$$
Therefore the likelihood ratio multiplied by minus two is
\begin{align*}
R(\alpha) &= -2 \log d\frac{\mathcal L_u(E,\Delta \vert \alpha)}{\mathcal L_u(E,\Delta \vert \hat{\alpha})}  \\
&= -2 \left[ \log\left( \frac{\alpha}{\hat{\alpha}}\right)\sum_{j=1}^{N_u} \delta_{[n-N_u+j]} -(\alpha -\hat{\alpha}) \sum_{j=1}^{N_u} \log e_{j,u} \right]\\
&= 2k \frac{\sum_{j=1}^{N_u} \delta_{[n-N_u+j]}}{k} \left[ \frac{\alpha}{\hat{\alpha}} - 1 - \log \left( \frac{\alpha}{\hat{\alpha}}\right) \right]
\end{align*}
The following corollary immediately follows from the fact that $\frac{\sum_{j=1}^{N_u} \delta_{[n-N_u+j]}}{k}\to p$ and $\sqrt{k} (\alpha -\hat{\alpha}) \xrightarrow{\mathcal{D}} \mathcal{N}(0,\frac{p}{\alpha^2})$ in \cite{amer}.

\begin{corollary}
    Assuming (\ref{eq1.5}) holds and under regularity condition in theorem 2.1 in \cite{amer}, then $R(\alpha)$ has an asymptotic $\chi_1^2$ distribution.
\end{corollary}
Based on the abode corollary, a $100(1-\theta)\%$ confidence interval for $\alpha$ is
$$
I_{ML}(\theta):=\lbrace \alpha : R(\alpha) \leqslant q_\theta \rbrace,
$$
where $q_\theta$ is the $\theta$-level critical point of $\chi_1^2$.

In Bayesian framework, this prior distribution serves as a foundation for incorporating our subjective beliefs into the statistical analysis. It encapsulates our knowledge or assumptions about the parameter of interest, allowing us to update our beliefs based on observed data. In our setting, we provide the unknown tail index $\alpha$ with a Jeffrey prior density $\pi(\alpha)$. Then, using Bayes theorem, we obtain the posterior density $\pi(\alpha \vert E,\Delta)=\mathcal{L}_u(E,\Delta \vert \alpha)\pi(\alpha)\slash \int_\Lambda {\mathcal{L}_u(E,\Delta \vert \alpha)\pi(\alpha)d\alpha}$ of $\alpha$, where $\Lambda$ is the support of the distribution of $\alpha$. The posterior distribution of $\alpha$ is proportional to the product of the partial likelihood \eqref{partlik} and the prior, namely: $\pi(\alpha \vert E,\Delta)\propto \mathcal{L}_u(E,\Delta \vert \alpha)\pi(\alpha)$. Choosing the prior density is a central issue in Bayesian estimation.  When  information available for prior elicitation is minimal, one can use objective (or non-informative) priors, such as Jeffrey's prior (see, e.g, \cite{Jeff:61}).

This prior is proportional to the square root of Fisher's information. If $\mathcal{L}_u(E,\Delta \vert \alpha)$ is the likelihood of a single observation $(E, \Delta)$, Jeffrey's prior can be written as:
\begin{eqnarray}\label{jprior}
\pi (\alpha)&\propto &\left[ -\mathbb{E}\left(\dfrac{\partial ^2}{\partial \alpha ^2} \log \mathcal{L}_u(E,\Delta \vert \alpha)\right) \right]^{\frac{1}{2}} \nonumber\\
&\propto & \frac{1}{\alpha}.
\end{eqnarray}
Using \eqref{partlik} and \eqref{jprior}, the posterior density of $\alpha$ based on Jeffrey's prior is given by:
\begin{eqnarray*}
\pi (\alpha \vert E,\Delta) \propto \alpha^{(\sum_{j=1}^{N_u}{\delta_{[n-N_u+j]}})-1} \exp \left[ -\alpha \sum_{j=1}^{N_u}{\log e_{j,u}}\right],
\end{eqnarray*}
which coincides (up to some normalizing constants) with the probability density function of the distribution $Gamma(\sum_{j=1}^{N_u}{\delta_{[n-N_u+j]}},\sum_{j=1}^{N_u}{\log e_{j,u}})$. Based on this posterior, we construct two classical Bayesian estimators of $\alpha$, namely the mean posterior estimator (MPE) and the maximum posterior estimator (MAP).

Letting $u=Z_{n,n-k}$, we obtain the following formal Bayesian estimators of the tail index $\alpha$:

\begin{equation}\label{MPEJeff}
\hat{\alpha}_{MPE}^{(J)}:=\dfrac{\sum_{i=1}^{k}{\delta_{[n-i+1]}}}{\sum_{i=1}^{k}{\log \left(\dfrac{Z_{n,n-i+1}}{Z_{n,n-k}}\right) }},
\end{equation}

and

\begin{equation}\label{MAPJeff}
\hat{\alpha}_{MAP}^{(J)}:=\dfrac{\sum_{i=1}^{k}{\delta_{[n-i+1]}} -1}{\sum_{i=1}^{k}{\log \left( \dfrac{Z_{n,n-i+1}}{Z_{n,n-k}}\right) }}.
\end{equation}

A $100(1-\theta)\%$ highest posterior density interval (HPDI)  for $\alpha$ is defined as:
$$
I_{B}(\theta)= \min_{|a-b|} \lbrace (a,b)\in \mathbb{R}^+ : \mathbb{P}\left(\alpha \notin (a,b) \right) \leqslant \theta\rbrace,
$$
where $\mathbb{P}$ is the probability function associated to the $Gamma\left(\sum_{i=1}^{k}\delta_{[n-i+1]},\sum_{i=1}^{k}{\log \left( \dfrac{Z_{n,n-i+1}}{Z_{n,n-k}}\right) }\right)$ distribution.

Instead of using the complete likelihood, composite likelihood is an alternative constructed when partial likelihoods is considered on subsets of data or omitting certain components from the full likelihood. The composite likelihood is more appropriate framework than the partial likelihood in (\ref{partlik}), and provides an estimation model of parameters without relying on the entire likelihood function. Based on the sample $(E,\Delta)=(e_{j,u},\delta_{[n-N_u+j]})_{j=1,\ldots,N_u}$, the composite likelihood of $\alpha$ is given by:
\begin{eqnarray}\label{complik}
\mathcal L_u^{CL}(E,\Delta \vert \alpha)&=&\prod_{j=1}^{N_u}{\left[ \alpha {e_{j,u}}^{-(\alpha +1)}\right] ^{N_u w_i\delta_{[n-N_u+j]}}\left[ {e_{j,u}}^{-\alpha}\right] ^{N_u w_i(1-\delta_{[n-N_u+j])}}} \nonumber \\
&=& \alpha ^{N_u\sum_{j=1}^{N_u}{w_i\delta_{[n-N_u+j]}}}\left( \prod_{j=1}^{N_u}{e_{j,u}}^{N_u w_i}\right)^{-\alpha}\left( \prod_{j=1}^{N_u}{{e_{j,u}}^{-N_u w_i\delta_{[n-N_u+j]}}}\right),
\end{eqnarray}
where $w=(w_1,\ldots,w_{N_u})$ is a vector of non-negative weights with summation equal to one. Under a Jeffrey prior, the posterior distribution of $\alpha$ can be derived form the proportional relation
\begin{eqnarray}\label{comppost}
\pi_u^{CL}(\alpha \vert E,\Delta  )&\propto& \alpha ^{N_u(\sum_{j=1}^{N_u}{w_i\delta_{[n-N_u+j]}})-1}\left( \prod_{j=1}^{N_u}{e_{j,u}}^{N_u w_i}\right)^{-\alpha}, \nonumber \\
&\propto& \alpha^{N_u(\sum_{j=1}^{N_u}{w_i\delta_{[n-N_u+j]}})-1}\exp \left( -\alpha N_u\sum_{j=1}^{N_u}{ w_i \log e_{j,u}}\right), \nonumber \\
&\sim& Gamma\left(N_u\sum_{j=1}^{N_u}{w_i\delta_{[n-N_u+j]}}, N_u\sum_{j=1}^{N_u}{ w_i \log e_{j,u}}\ \right)
\end{eqnarray}
Letting $u=Z_{n,n-k}$, we obtain the following Bayesian-CL estimators of the tail index $\alpha$:

\begin{equation}\label{MPECL}
\hat{\alpha}_{MPE}^{(CL)}:=\dfrac{\sum_{i=1}^{k}{w_i\delta_{[n-i+1]}}}{\sum_{i=1}^{k}{w_i\log \left(\dfrac{Z_{n,n-i+1}}{Z_{n,n-k}}\right) }},
\end{equation}

and

\begin{equation}\label{MAPCL}
\hat{\alpha}_{MAP}^{(CL)}:=\dfrac{\sum_{i=1}^{k}{w_i\delta_{[n-i+1]}} -\frac{1}{k}}{\sum_{i=1}^{k}{w_i\log \left( \dfrac{Z_{n,n-i+1}}{Z_{n,n-k}}\right) }}.
\end{equation}

Based on (\ref{comppost}), a $100(1-\theta)\%$ highest posterior density interval (HPDI)  for $\alpha$ is defined as:
$$
I_{BCL}(\theta)= \min_{|a-b|} \lbrace (a,b)\in \mathbb{R}^+ : \mathbb{P}_{CL}\left(\alpha \notin (a,b) \right) \leqslant \theta\rbrace,
$$
where $\mathbb{P}_{CL}$ is the probability function associated to the $Gamma\left(k\sum_{i=1}^{k}w_i\delta_{[n-i+1]},k\sum_{i=1}^{k}{w_i\log \left( \dfrac{Z_{n,n-i+1}}{Z_{n,n-k}}\right) }\right)$ distribution.

\cite{lu02} define a distance function $\mathbb{Q}(w)=\frac{1}{k}\sum_{i=1}^k  \log(k w_i)$, which is the measure of distance between $w$ and uniform distribution $(w_i)_{i=\overline{1,k}}=1/k$. Our estimation strategy consists in choosing the weights $w_i$ to minimize the distance $(2k)^{-1}\mathbb{L}(\alpha)=\mathbb{Q}(w)$ according to the constraints
$$
\sum_{i=1}^k w_i =1 \quad \mbox{ and } \quad k\sum_{i=1}^k\dfrac{ w_i \hat{p}_k}{V_i} =\alpha,
$$
where $\hat{p}_k=\frac{1}{k}\delta_{[n-i+1]}$ and $V_i=\log\left(Z_{n,n-i+1}/Z_{n,n-k}\right)$ for $i=1,\ldots, k$. The log-posterior ratio is defined as $\ell(\alpha)=-2\log\left(k^k \mathbb{L}(\alpha) \right)$ and by the Lagrange's method, we obtain
$$
w_i= \frac{1}{k}\left[1+\lambda \left( \frac{k\hat{p}_k}{V_i}-\alpha\right)\right]^{-1}, \quad \mbox{ and } \quad \ell(\alpha)=2\sum_{i=1}^k{ \log \left( 1+ \lambda \left( \frac{k\hat{p}_k}{V_i}-\alpha\right) \right)},
$$
where  $\lambda$ is solution of the equation
\begin{equation}
\sum_{i=1}^k{\left( \frac{k\hat{p}_k}{V_i}-\alpha\right) \left[1+\lambda \left( \frac{k\hat{p}_k}{V_i}-\alpha\right)\right]^{-1}} =0
\label{eq.lambda}
\end{equation}
\begin{corollary}\label{cor01}
     Assuming (\ref{eq1.5}) holds and under regularity condition in theorem 2.1 in \cite{amer}, then $\ell(\alpha)$ has an asymptotic $\chi_1^2$ distribution.
\end{corollary}
based on the above corollary, a $100(1-\theta)\%$ confidence interval for the tail index $\alpha$ is
$$
I_{EL}(\theta):=\lbrace \alpha : \ell(\alpha) \leqslant q_\theta \rbrace,
$$
where $q_\theta$ is the $\theta$-level critical point of $\chi_1^2$.

\section{Simulation study}\label{sec4}

In this section, we assess, via simulations, finite-sample performance of confidence regions proposed in section \ref{sec3} (The consistency of confidence interval $I(\theta)$ defined in section \ref{sec2} will be performed through real-data study). 
\vspace{.2cm}

\noindent \textit{Study design.} Our simulation design is as follows: Let $X$ and $Y$ be independent random variables with cdf $F$ and $G$ respectively satisfying (\ref{eq.1}), with tail index $\alpha>0$ and $\beta>0$ respectively. As $u \rightarrow \infty$, $p_u:=\mathbb P(\delta=1|Z>u)$ tends to $\frac{\alpha}{\alpha+\beta}$. Thus, for a given value of $\alpha$, a relevant choice of $\beta$ allows us to generate data with an approximate proportion of non-censored $X_i$s among all $Z_i >u$. The simulation procedure is replicated as follows:

\begin{enumerate}
\item \textit{Generate a sample of  $n$  independent copies of $(Z, \delta)$, where $Z=\min(X,Y)$, $\delta= \mathds{1}_{\left\lbrace X \leq Y \right\rbrace}$. Given $\alpha$, and varying the value of $\beta$ to allow censoring percentage in the right tail of $X$ to be approximately $5\%, 10\%, 30\%$ and $50\%$.}
\item \textit{For each $\beta$, we compute the proposed confidence regions $I_{ML}$, $I_{B}$, $I_{BCL}$ and $I_{EL}$, for a given level $\theta$, and by incrementing the fraction level $\frac{k}{n}$ from $\frac{k_{min}}{n}$ to $\frac{k_{max}}{n}$. }
\item \textit{Steps 1-2 are repeated $m$ times, so that we obtain $m$ realisations of each $I_{\bullet}(k)$ (where $I_{\bullet}(k)$ is any of the confidence regions calculated in step 2), for each  $\beta$.}
\item \textit{For each $\beta$, we return coverage probabilities as the proportions of
the N intervals which contain the true value of $\alpha$. Also, we indicate the average proportion of uncensored data among the observations $Z_i > Z_{n,n-k}$ and the average length of each $I_{\bullet}(k)$ (for $k=k_{min},\ldots,k_{max}$) over the $m$ replications.}
\end{enumerate}
To implement the procedure cited above, we consider the following two simulation settings:
\begin{itemize}
\item[] \textbf{Example 01:} $X$ follows a standard Generalized Pareto distributions (GPD) and $Y$ follows a Fr\'echet distributed, with cdf $F_{\mu , \sigma ,\alpha}(x)=1-(1+\frac{x}{\alpha})^{-\alpha}$ and $G_{\beta}(x)=\exp (-x^{-\beta })$ respectively,
\item[] \textbf{Example 02:} $X$ and $Y$ are both distributed as standard log-logistic random variables with cdf $F_{\alpha}(x)=\frac{1}{1+x^{-\alpha}}$ and $G_{\beta}(x)=\frac{1}{1+x^{-\beta}}$ respectively (the log-logistic distribution is commonly used in health science to model survival data).
\end{itemize}
\vspace{.2cm}

\begin{figure}[H]
\begin{center}
\includegraphics[keepaspectratio=true,scale=0.42]{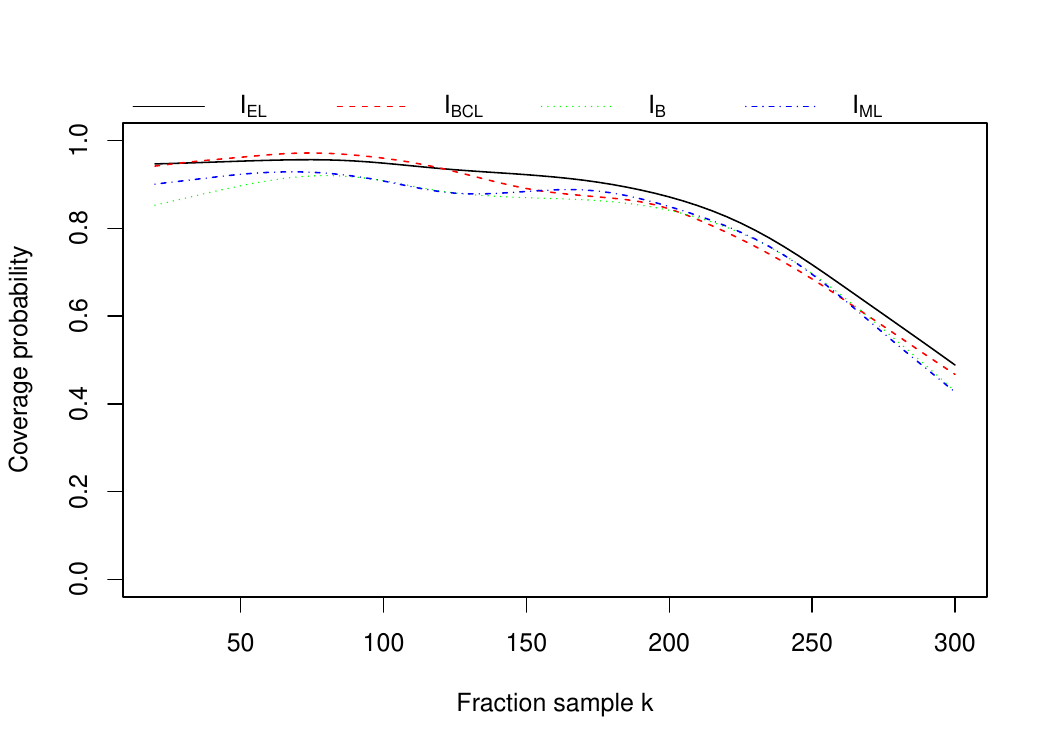}\hfill 
\includegraphics[scale=0.42]{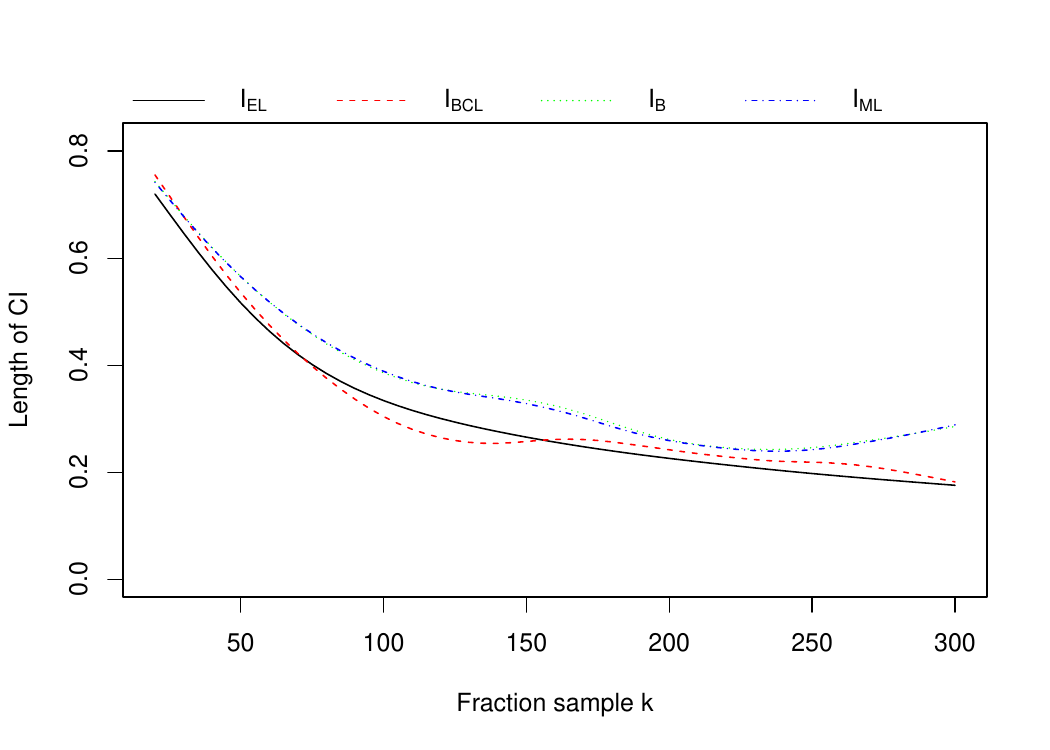}\\
\includegraphics[keepaspectratio=true,scale=0.42]{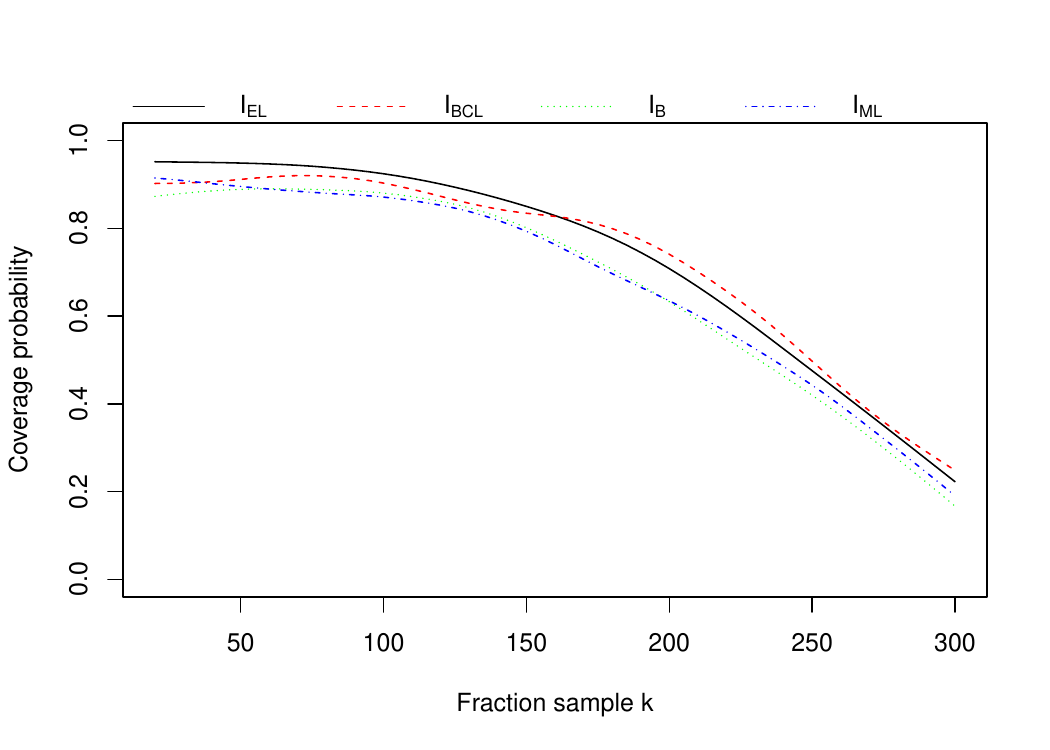}\hfill
\includegraphics[scale=0.42]{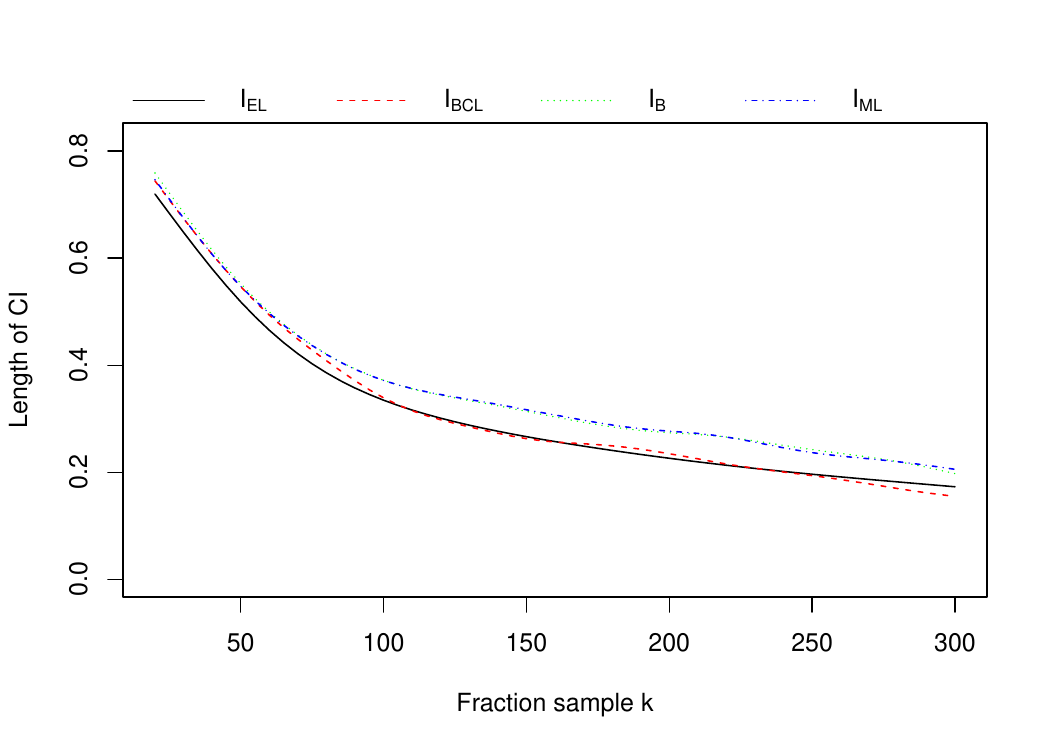}\\
\includegraphics[keepaspectratio=true,scale=0.42]{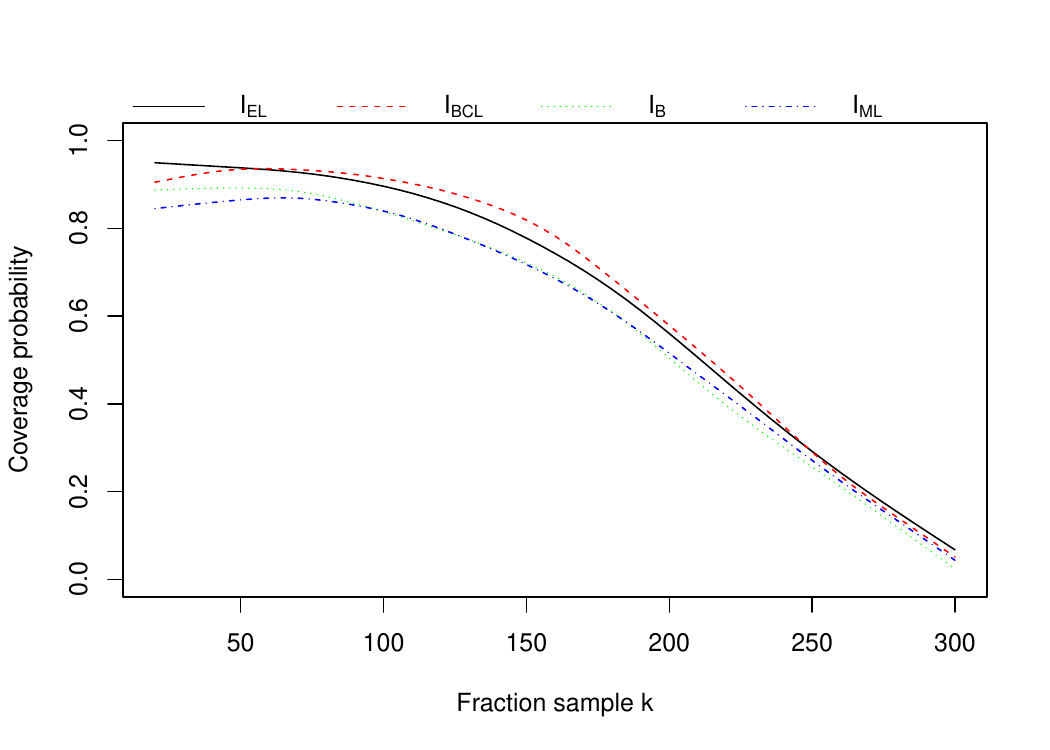}\hfill
\includegraphics[scale=0.42]{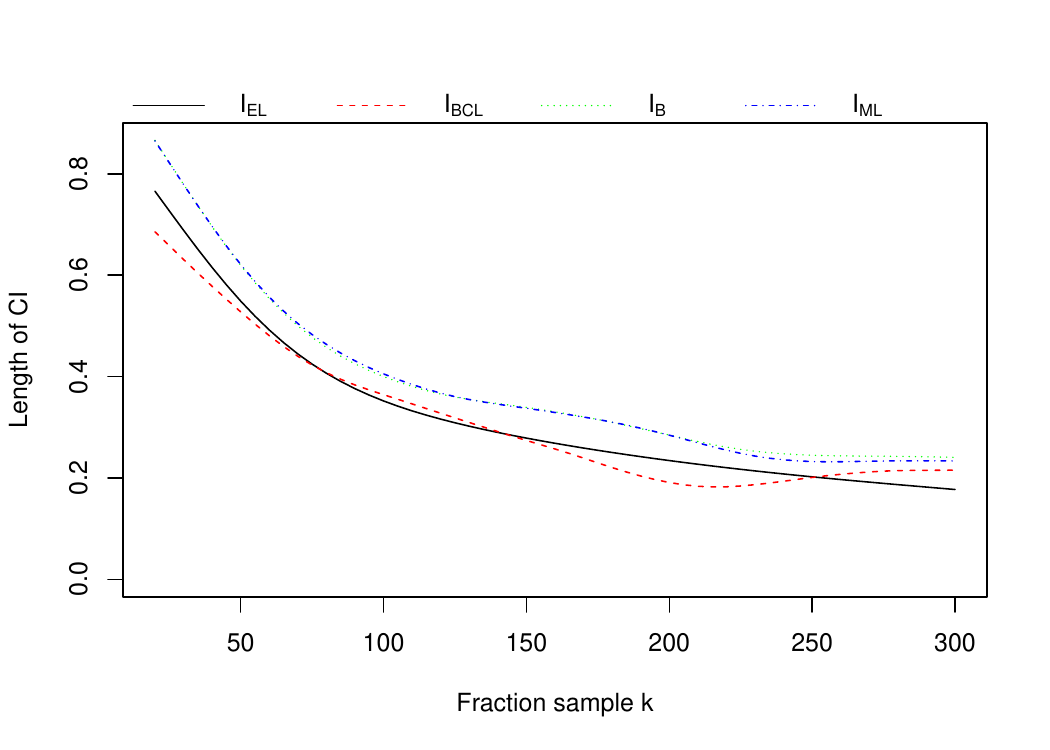}\\
\includegraphics[keepaspectratio=true,scale=0.42]{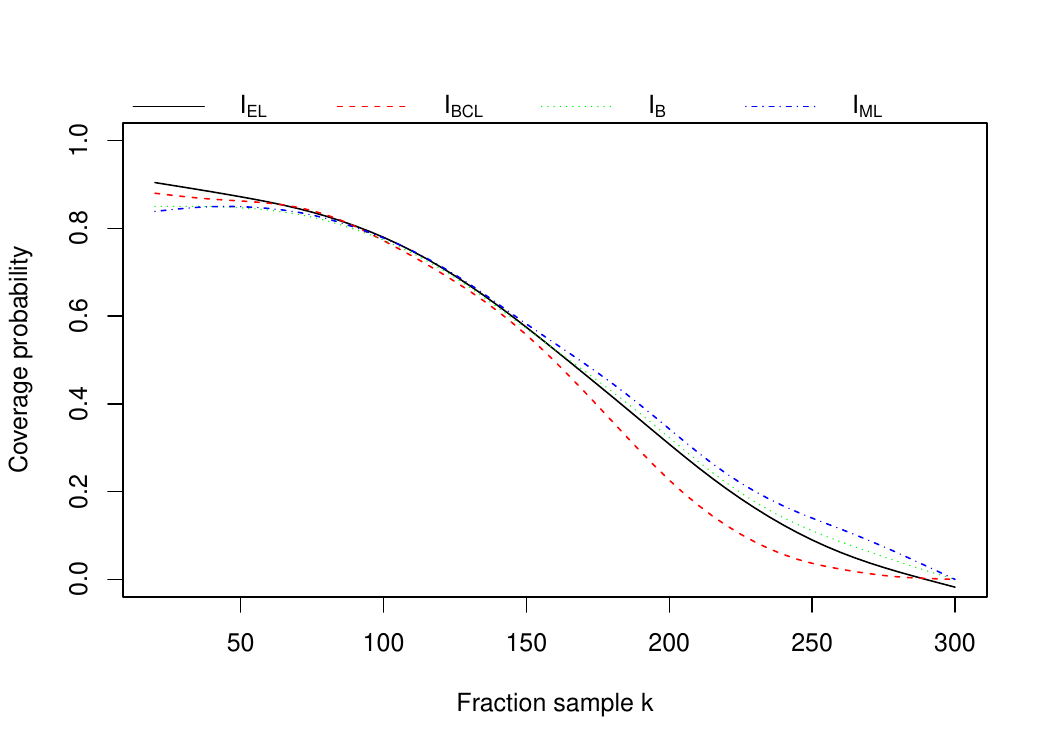}\hfill
\includegraphics[scale=0.42]{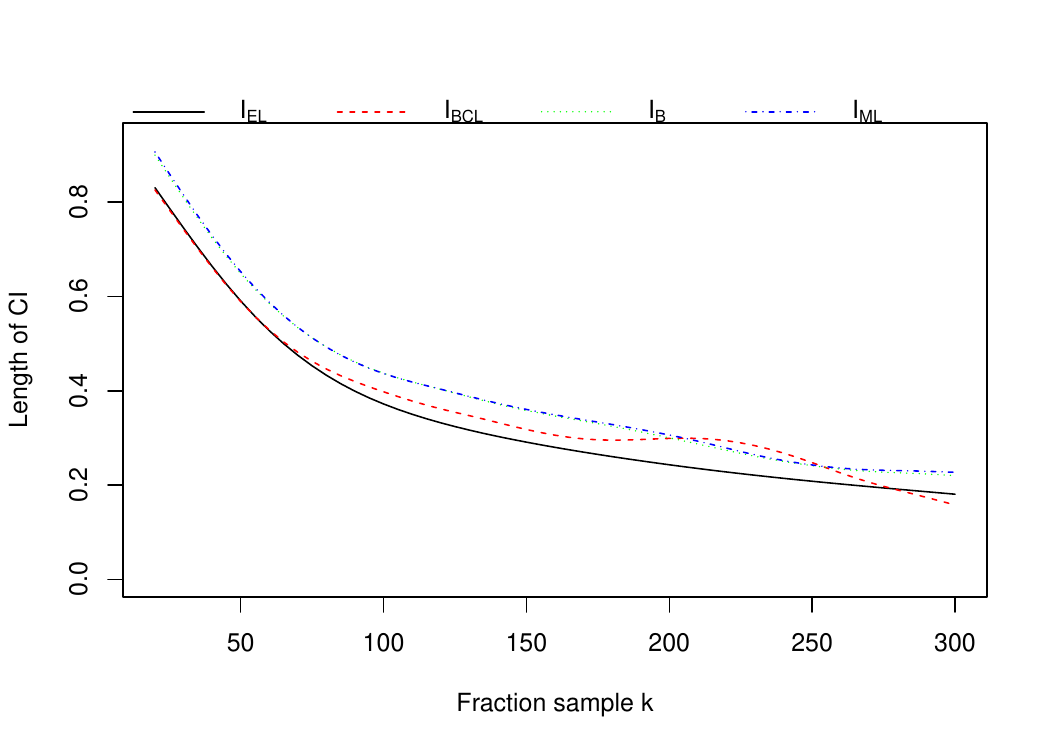}
\end{center}
\caption{Coverage probability (left) and average length of confidence region  (right) of $I_{EL}$ (black line), $I_{BCL}$ (red dashed), $I_{B}$ (green dotted),  $I_{ML}$ (blue dotdash) for GPD distribution with $\alpha=1.25$ censored by Fr\'echet distribution with $\beta=\frac{5}{100}$ (top row), $\beta=\frac{1}{10}$ (second row), $\beta=\frac{1}{2}$ (third row) and $\beta=1$ (bottom row).}
\label{figss1}
\end{figure}

\noindent \textit{Results for Example 01.} We consider $\alpha=1.25$ and $\beta=\frac{1}{100}, \frac{1}{10}, \frac{1}{2}$ and $\beta=1$. Simulations are conducted using the statistical software \texttt{R} \cite{softR}. Results are provided for a sample of size $n=1000$ and $m=1000$ simulated samples, with $k$ in the range of $k_{min}=30$ to $k_{max}=300$. We plot  the coverage probability and the average length of each confidence region versus $k$ (Figure \ref{figss1}). In Table \ref{tab:table1}, we report: the averaged (over the $m$ simulated samples) confidence regions length, the empirical coverage probability ($CP(k)$) of each $I_{\bullet}$ at the optimal fraction level $k_{opt}=\arg\min_{k}CP[I_{\bullet}(k)]$. For each $\hat{\alpha}_\ell$, and the averaged proportion $\bar{p}(k_{opt})=\frac{1}{k_{opt}}\sum_{i=1}^{k_{opt}}\delta _{[n-i+1]}$ of uncensored data among observations $i$ such that $Z_i >Z_{n,n-k_{opt}}$. For every $\beta$, the average value of $1-\bar{p}(k_{opt})$ is close to the target censoring proportion in the right tail of $X$ (namely, $0.05, 0.10, 0.30, 0.50$). From these results, it appears, as expected, that the coverage probabilities of all confidence regions decrease when censoring in the right tail of $X$ increases. Likewise, the confidence regions constructed under Bayesian composite likelihood $I_{EL}$ outperforms $I_{BCL}$, $I_{B}$ and $I_{ML}$ in terms of coverage probability and length of the confidence region, for almost every $k$ and  might be regarded as the best among all three regions. 

\begin{figure}[H]
\begin{center}
\includegraphics[keepaspectratio=true,scale=0.42]{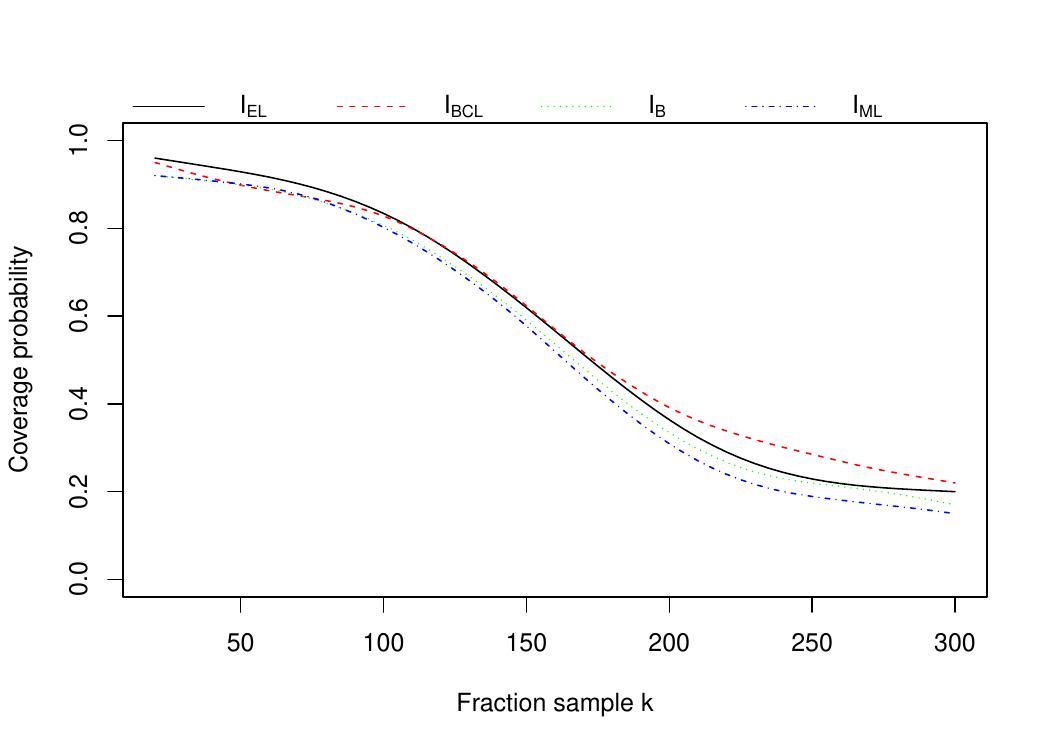}\hfill 
\includegraphics[scale=0.42]{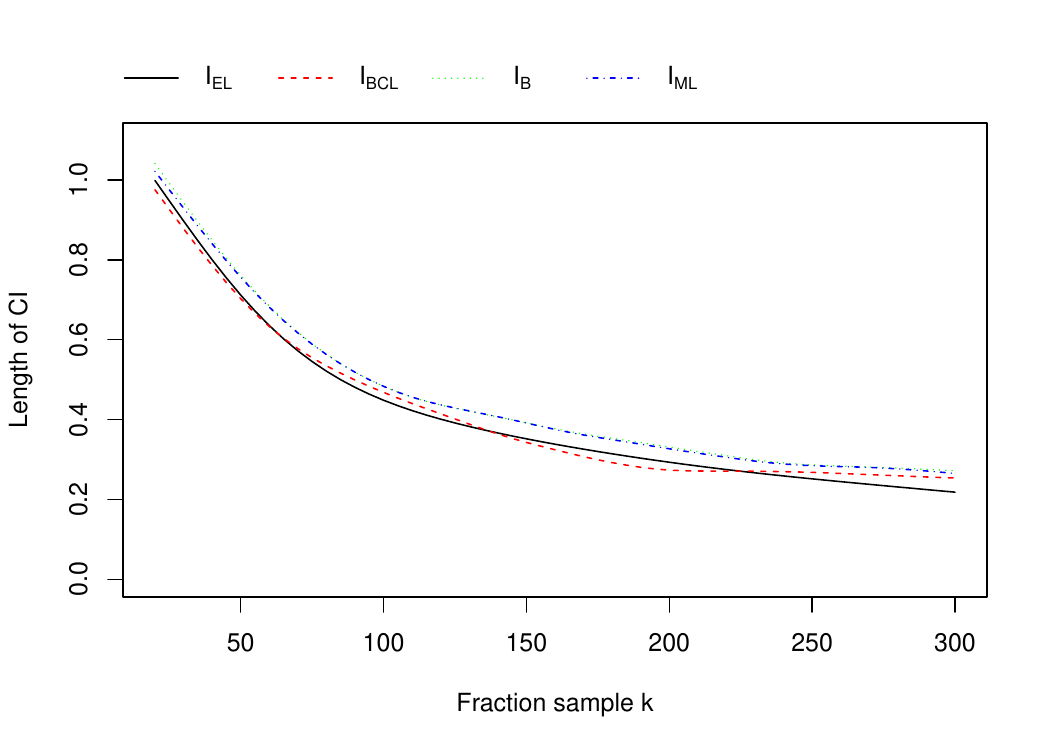}\\
\includegraphics[keepaspectratio=true,scale=0.42]{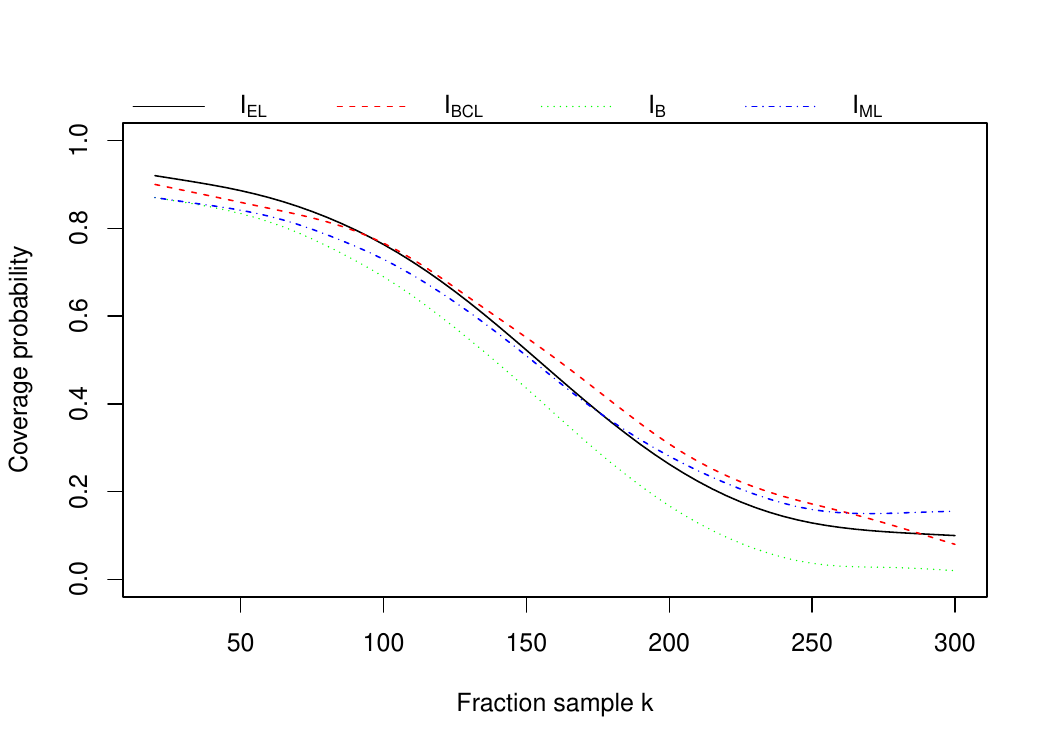}\hfill
\includegraphics[scale=0.42]{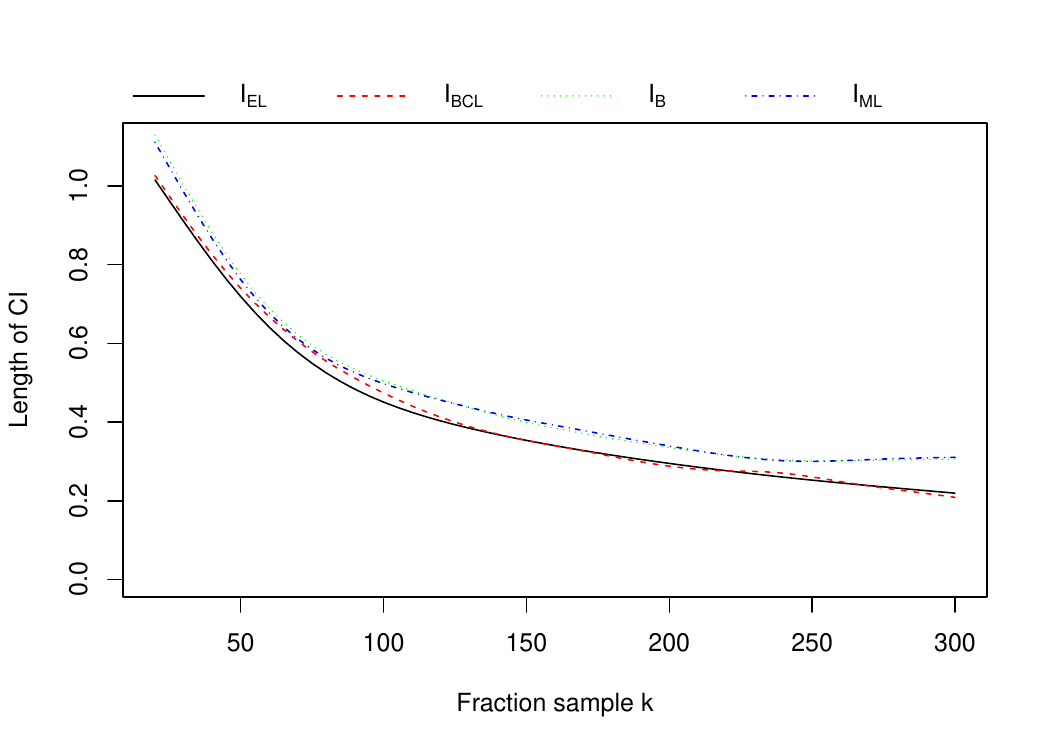}\\
\includegraphics[keepaspectratio=true,scale=0.42]{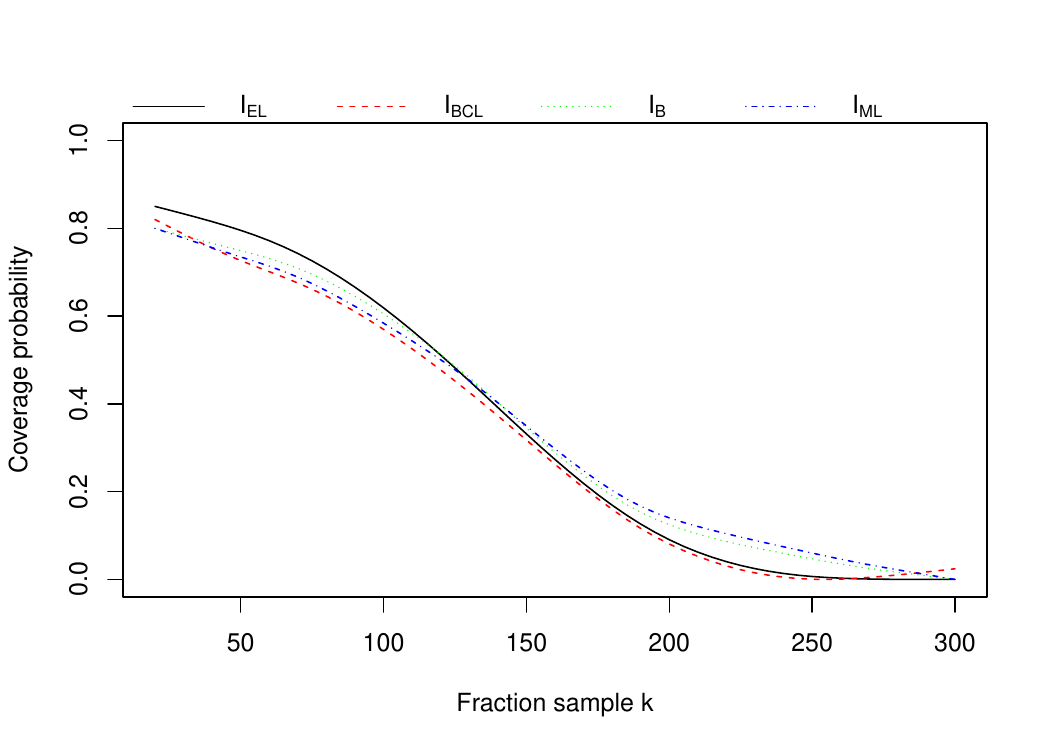}\hfill
\includegraphics[scale=0.42]{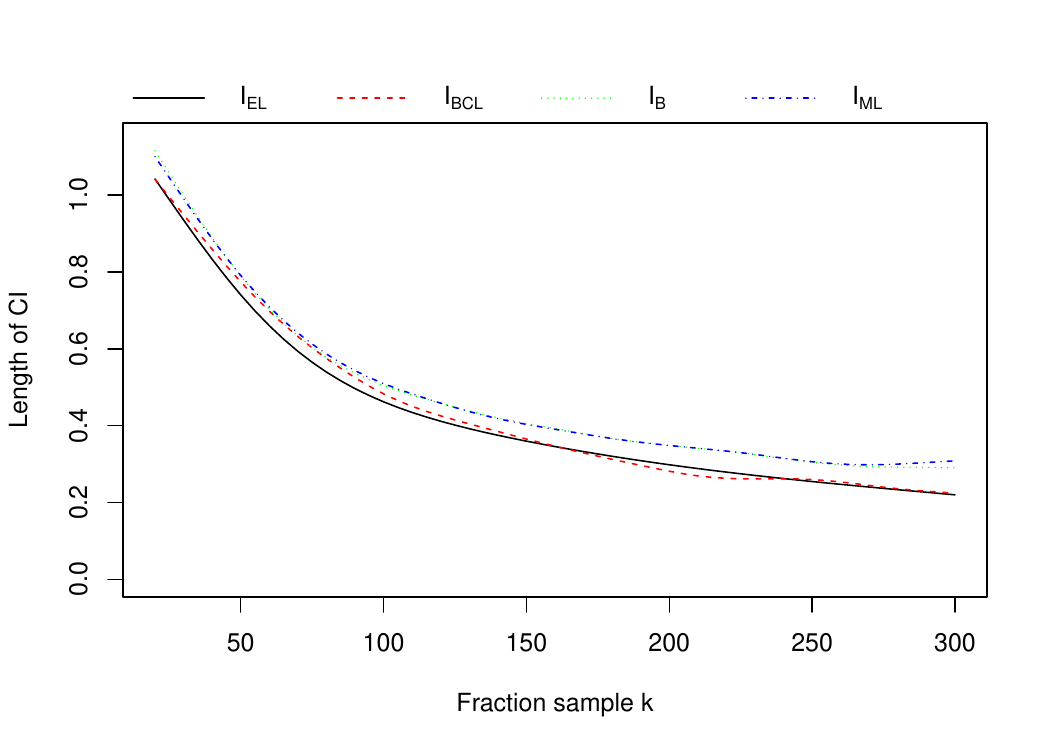}\\
\includegraphics[keepaspectratio=true,scale=0.42]{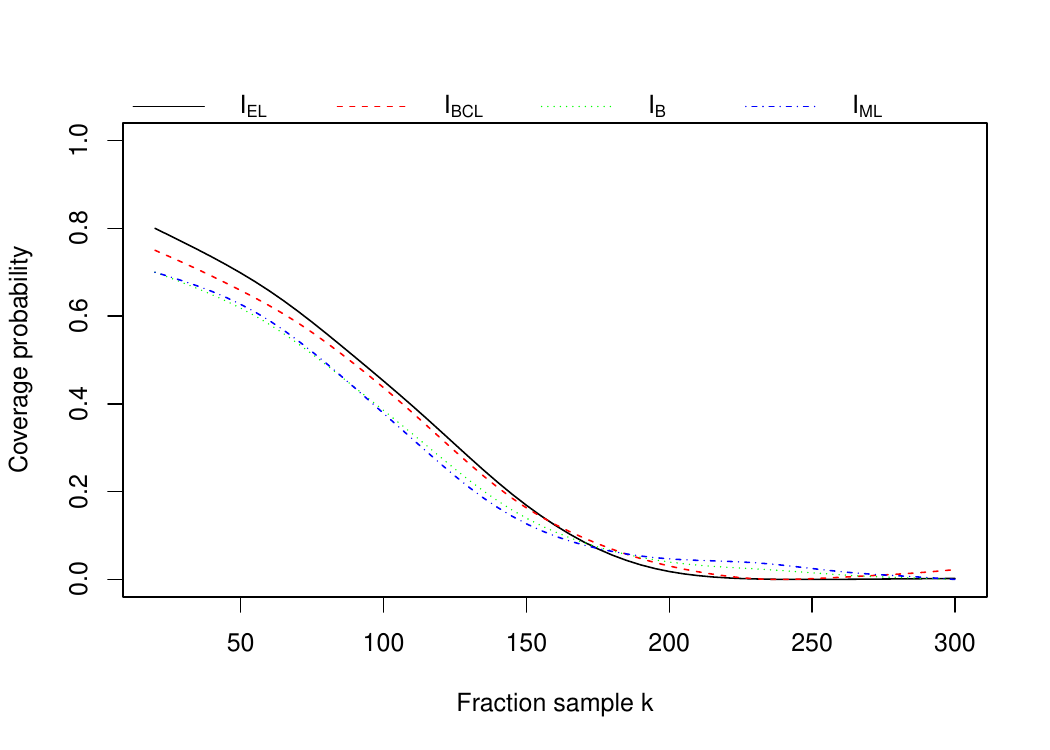}\hfill
\includegraphics[scale=0.42]{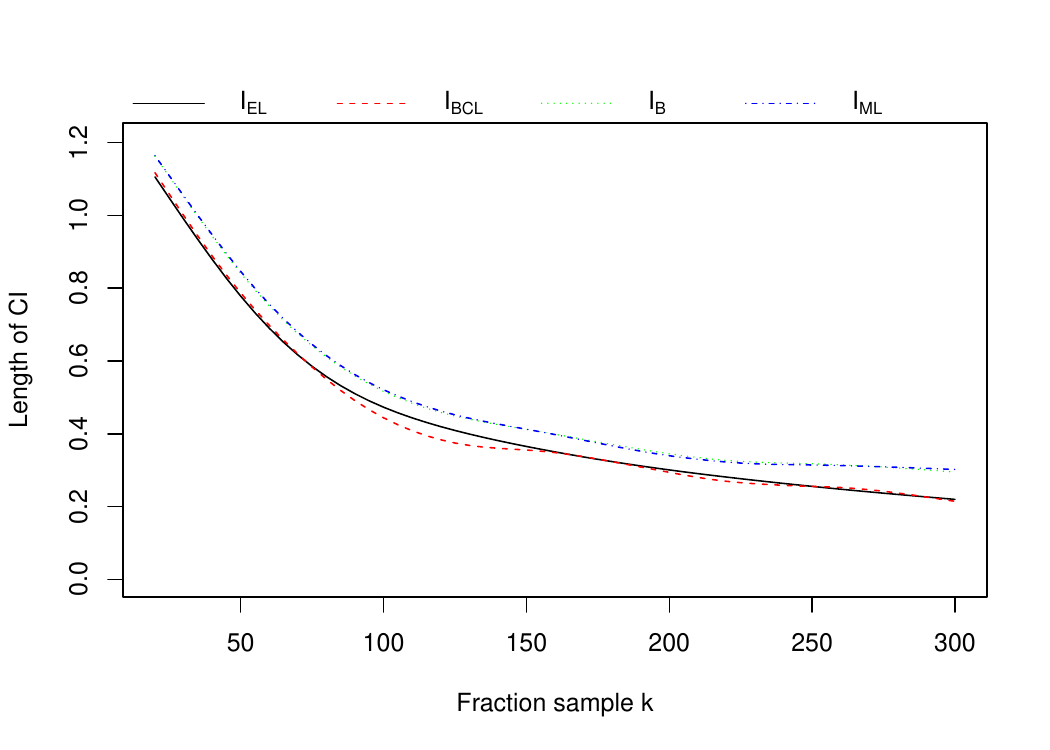}
\end{center}
\caption{Coverage probability (left) and average length of confidence region  (right) of $I_{EL}$ (black line), $I_{BCL}$ (red dashed), $I_{B}$ (green dotted),  $I_{ML}$ (blue dotdash) for Log-logistic distribution with $\alpha=1.25$ censored by Log-logistic distribution with $\beta=\frac{5}{100}$ (top row), $\beta=\frac{1}{10}$ (second row), $\beta=\frac{1}{2}$ (third row) and $\beta=1$ (bottom row).}
\label{figss2}
\end{figure}

\noindent \textit{Results for Example 02.} 
Letting $\alpha=0.85$ (tail index $\alpha<1$, coincides with significant thickness in the tail of the distribution) and $\beta=\frac{1}{100}, \frac{1}{10}, \frac{1}{2}$ and $\beta=1$, such that the expected proportion of uncensored in the tail of $X$ to be approximately $5\%, 10\%, 40\%$ and $55\%$. From Figure (\ref{figss2}),
the same conclusion can be drawn as in the previous example, but with more amplified effect of degradation of the coverage probabilities when the proportion of censorship increases.

\begin{table}
    \caption{Optimal results for GPD model with $\alpha =1.25$, censored by a Fr\'echet distribution for each $\beta$, the estimator with highest  empirical coverage is indicated in bold.}
	\centering
    \begin{tabular}{cp{01cm}p{0.8cm}p{1.5cm}p{2cm}p{2cm}}
    \hline
    \hline
          &       &  $k_{opt}$  &  $\bar{p}(k_{opt})$  &  CP  &  $\Bar{l}_{CI}$  \\
    \hline
    \multicolumn{1}{c}{\multirow{4}{*}{\textbf{$\beta=\frac{5}{100}$}}} & $I_{EL}$    & 49 & 0.9591 & 0.9571 & 0.5328 \\
    \multicolumn{1}{c}{} & $I_{BCL}$    & 87 & 0.9540 & \textbf{0.9628} & 0.3274 \\
    \multicolumn{1}{c}{} & $I_{B}$    & 73 & 0.9589 & 0.9264 & 0.4378 \\
    \multicolumn{1}{c}{} & $I_{ML}$    & 69 & 0.9420 & 1.1841 & 0.0148 \\
    \hline
    \multicolumn{1}{c}{\multirow{4}{*}{\textbf{$\beta=\frac{1}{10}$}}} & $I_{EL}$    & 37 & 0.9189 & \textbf{0.9405} & 0.6792 \\
    \multicolumn{1}{c}{} & $I_{BCL}$    & 81 & 0.9135 & 0.9264 & 0.4189 \\
    \multicolumn{1}{c}{} & $I_{B}$    & 35 & 0.9142 & 0.8926 & 0.6902 \\
    \multicolumn{1}{c}{} & $I_{ML}$    & 78 & 0.9102 & 0.8764 & 0.4275 \\
    \hline
    \multicolumn{1}{c}{\multirow{4}{*}{\textbf{$\beta=\frac{1}{2}$}}} & $I_{EL}$    & 36 & 0.6944 & \textbf{0.9387} & 0.6981 \\
    \multicolumn{1}{c}{} & $I_{BCL}$    & 53 & 0.7169 & 0.9268 &  0.5134\\
    \multicolumn{1}{c}{} & $I_{B}$    & 48 & 0.7083 & 0.8722 & 0.4376 \\
    \multicolumn{1}{c}{} & $I_{ML}$    & 74 & 0.7162 & 0.8653 & 0.4376 \\
    \hline
    \multicolumn{1}{c}{\multirow{4}{*}{\textbf{$\beta=1$}}} & $I_{EL}$    & 34 & 0.5294 & \textbf{0.9021} & 0.8027 \\
    \multicolumn{1}{c}{} & $I_{BCL}$    & 34 & 0.5294 & 0.8679 & 0.8035 \\
    \multicolumn{1}{c}{} & $I_B$    & 58 & 0.5172 & 0.8136 & 0.5943 \\
    \multicolumn{1}{c}{} & $I_{ML}$    & 58 & 0.5172 & 0.8094 & 0.5962 \\
    \hline
    \hline
    \end{tabular}%
	\label{tab:table1}
\end{table}

\section{Real data application}\label{sec5}
\subsection{Global most costly natural disasters (1900--2024)}
In this subsection, we explore the performance of the the confidence interval proposed in section (\ref{sec2}), through a complete dataset of most expensive disaster in human history expressed in billion of dollars. The data are extracted from the website Our world in data (\href{https://ourworldindata.org/grapher/economic-damage-from-natural-disasters}{Economic damage by natural disaster type} - see for instance \cite{rit:22}). The data represent a set of 122 largest costly natural catastrophes that the humanity has record from 1900 to 2024, with costs over one billion dollars (earthquakes, storms, wildfires, droughts, floods, ...), in different regions of the globe.

\begin{figure}[H]
	\centering
\includegraphics[width=0.45\linewidth]{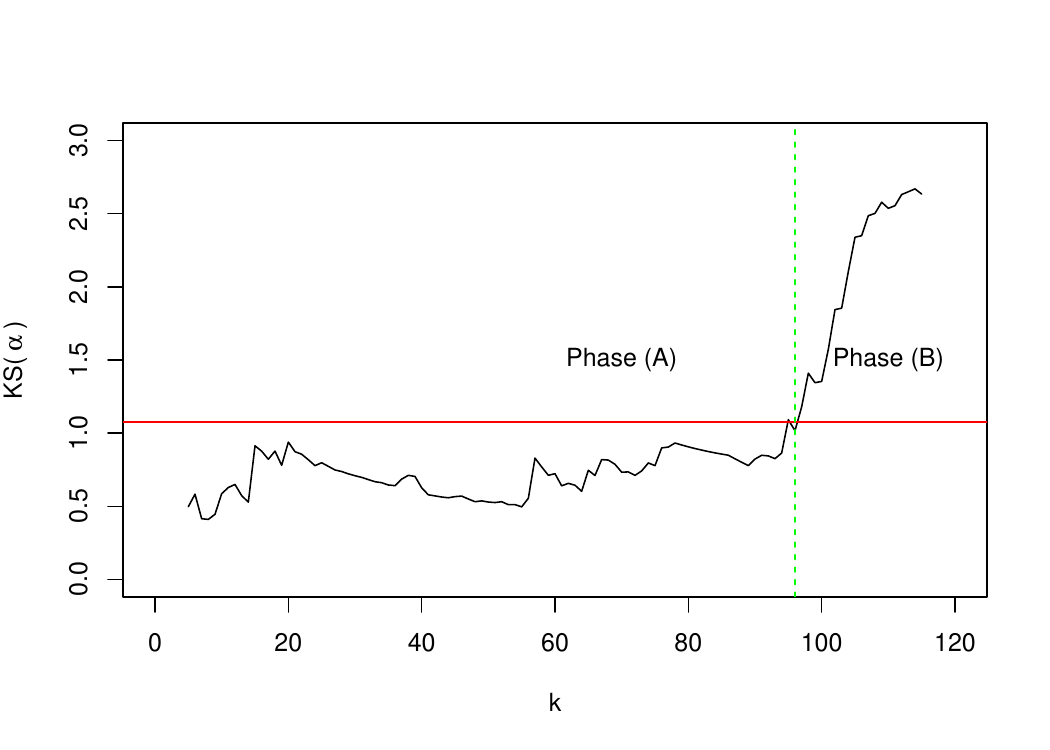}
\hfill\includegraphics[width=0.45\linewidth]{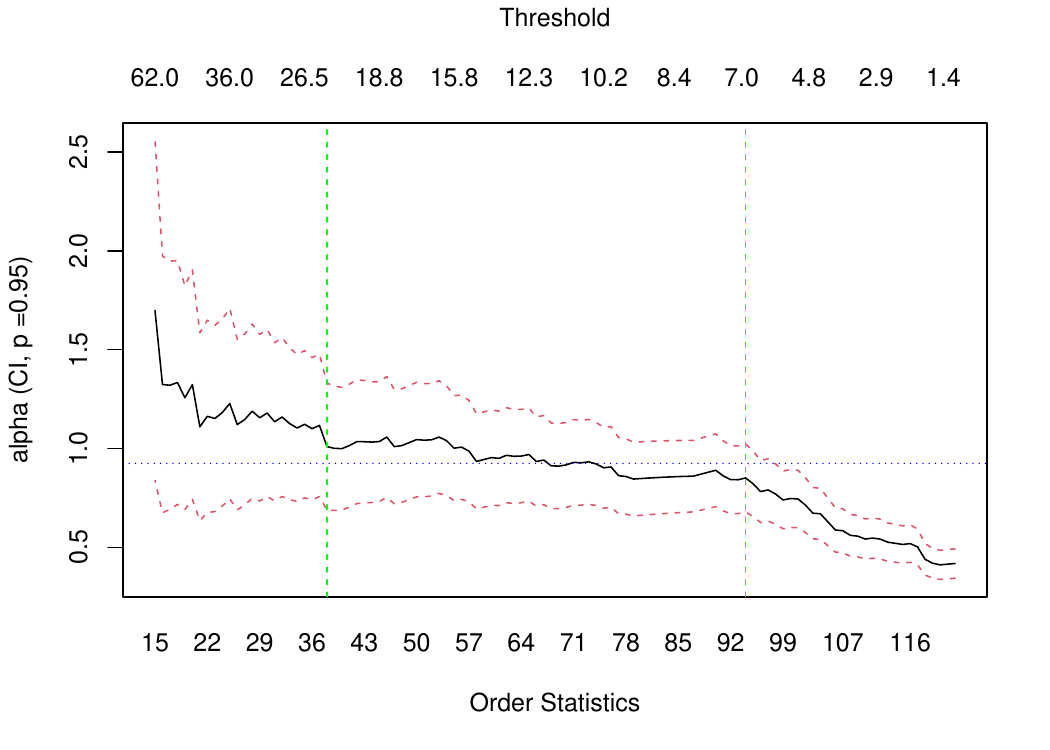}
\caption{Left : KS statistics  - Right : Hill estimate of the tail index $\alpha$ for natural disaster costs dataset}
	\label{figapp1}
\end{figure}

The presence of extreme values is a crucial aspect to be considered in real-data
analysis. Understanding these extremes can offer valuable insights and inform decision-making across various fields. An adapted Kolmogorov-Smirnov test for a heavy-tailed distribution was further
developed by \citep{Kon:08}, and confirmed if a heavy-tailed distribution is appropriate for a data set. For every $k = 1, 2, . . .$ the test
statistic is defined as:
$$
KS(k, \hat{\alpha}) = \sup_{r>1} \sqrt{ k} \vert 1 - G_k(r) - r^{\hat{\alpha}} \vert
$$
where $1-G_k(r) = \frac{1}{k} \sum_{i=1}^{n} \mathds{1}_{\{X_i > r X_{n,n-k}\}}$ and $\hat{\alpha}$ is the inverse of the Hill estimator in (\ref{eq.2}). 
\cite{Kon:08}, establish that the limiting critical value with level $0.95$ is $1.076$. This value is represented as the horizontal line on Figure (\ref{figapp1}). Koning and Peng (2008) explain that: \textit{"When the null hypothesis is true, test should  reject the null hypothesis for large values of $k$ since the critical values are obtained
by ignoring the bias"}. Figure (\ref{figapp1}) seems to be consistent with the ideal pattern described by \cite{Kon:08}, which suggests that the heavy-tailed
hypothesis holds for the natural disasters data. Indeed, for every $k\leqslant 95$
(phase (A)), the hypothesis of a heavy-tailed distribution is consistently accepted. When $k$ is larger than 95 (phase (B)), the hypothesis is consistently rejected.

We also plot the Hill estimate of the tail index of heavy-tailed data, to check which stable region of $k$ are suitable for the estimation of $\alpha$. From Figure (\ref{figapp1}), we deduce that the average value of the estimated tail index $\hat{\alpha}$ is around $0.9368$ (horizontal blue dotted line) for $k$ between $38$ and $95$.

\begin{figure}[H]
	\centering
\includegraphics[scale=0.45]{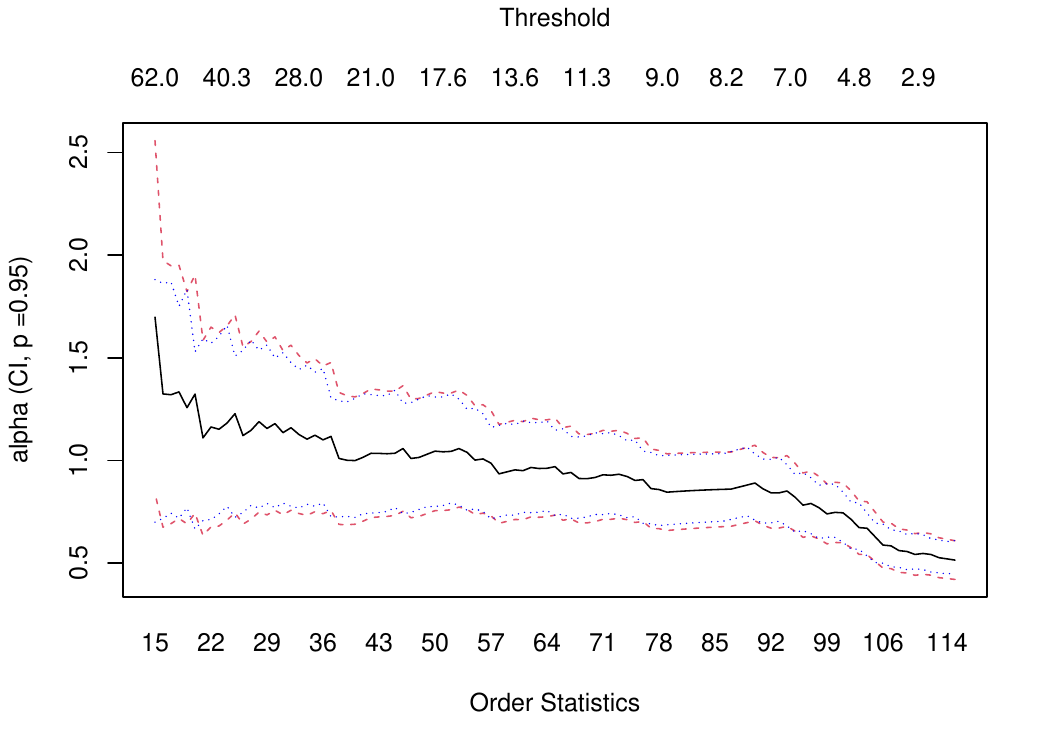}
\caption{Length of the approximate confidence regions of $\alpha$ with level $0.95$, for natural disaster costs data. Hill estimator (black line), Normal asymptotic approximation confidence interval (dashed red line) and $I(\theta)$ confidence interval (dotted blue line) are plotted against $k = 15, 16, \ldots , 120$.}
\label{figapp11}
\end{figure}

In Figure (\ref{figapp11}), we plotted the length of the approximate confidence regions $I(\theta)$  of $\alpha$ defined in (\ref{eq2.12}) and the asymptotic Gaussian confidence interval with level $\theta=0.95$ against the fraction $k = 15, 16, \ldots , 115$. It is worth highlighting that the $I(\theta)$ confidence interval consistently yields the smallest approximate interval lengths across the majority of cases according to different values of $k$.

\subsection{Australian AIDS survival data}

In this subsection, we investigate the performance of our methodology on a medical data. The dataset contains $n=2843$ patients diagnosed with AIDS in Australia before 1 July 1991. (see for more details \cite{VenRip:02}). The data are available in the \texttt R package "MASS" \cite{softR}.  $61.95\%$ of patients died. The other survival times are considered right-censored.

\cite{amer} analyzed this AIDS survival data by distinction of sex or age and explored the presence of extremes survival duration using the \cite{Kon:08} test. The hypothesis of a heavy-tailed distribution for $X$ is accepted by KS test when
$k$ ranges over the integer interval $\{20,\ldots , 218\}$. Our objective here is to illustrate the behavior of the proposed confidence regions for the tail index $\alpha$ with level $\theta=95\%$. To perform our application on the AIDS dataset, and since age and sex are important prognostic factors, we consider sub-populations of the whole data set. For illustrative purpose, we consider three sub-populations: male individuals (without distinction of age),  male individuals with age less than or equal to 35 years and male individuals with age greater than 35 years. These sub-populations are denoted by $\mathcal S_1, \mathcal S_2, \mathcal S_3$ respectively.

We calculate the length of the confidence regions $I_{EL}$, $I_{BCL}$, $I_{B}$ and $I_{ML}$ on $\mathcal S_1, \mathcal S_2, \mathcal S_3$. We also propose to plot the proportion $\hat{p}(k)=\frac{1}{k}\sum_{i=1}^n{\delta _i \mathds{1}_{\{Z_i > Z_{n,n-k}\}}}$ of non-censored observations in the $k$ largest $Z_i$'s against $k=30,31, \ldots, 500$.

\begin{figure}[h]
\begin{center}
\includegraphics[keepaspectratio=true,scale=0.42]{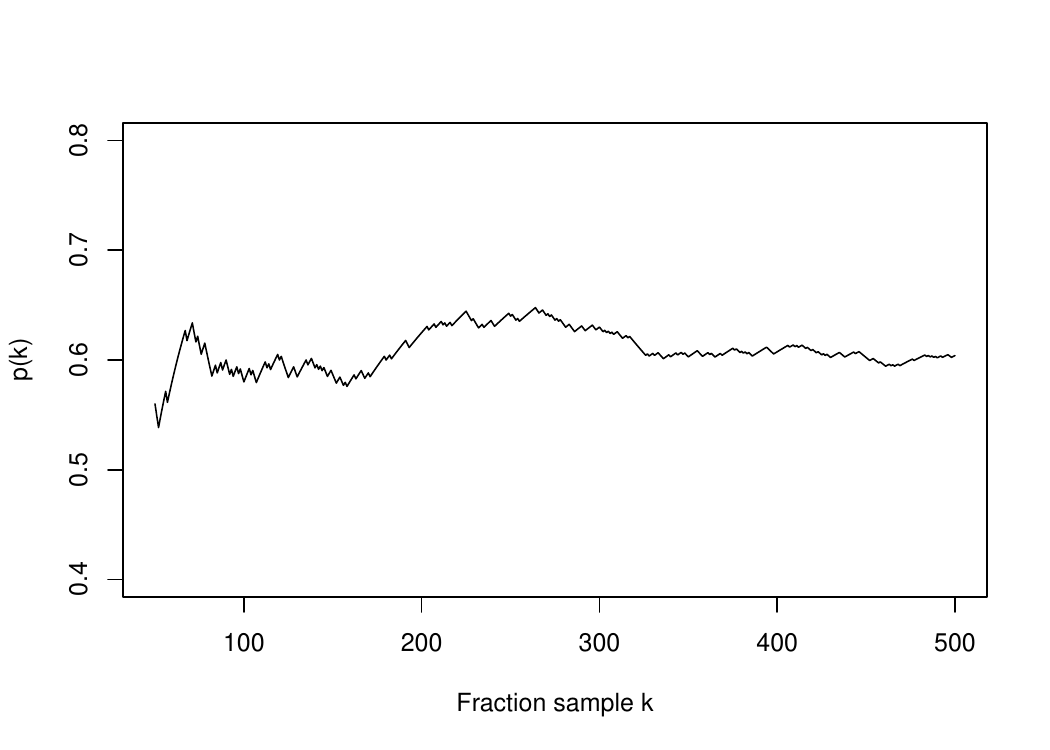}\hfill 
\includegraphics[scale=0.42]{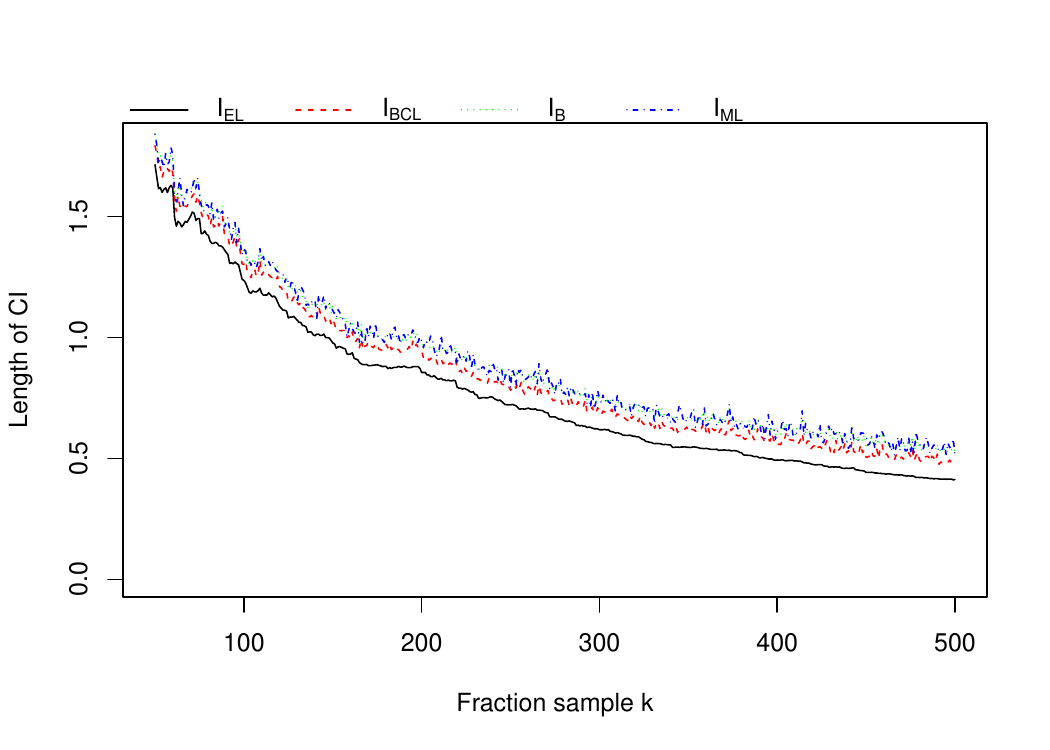}\\
\includegraphics[keepaspectratio=true,scale=0.42]{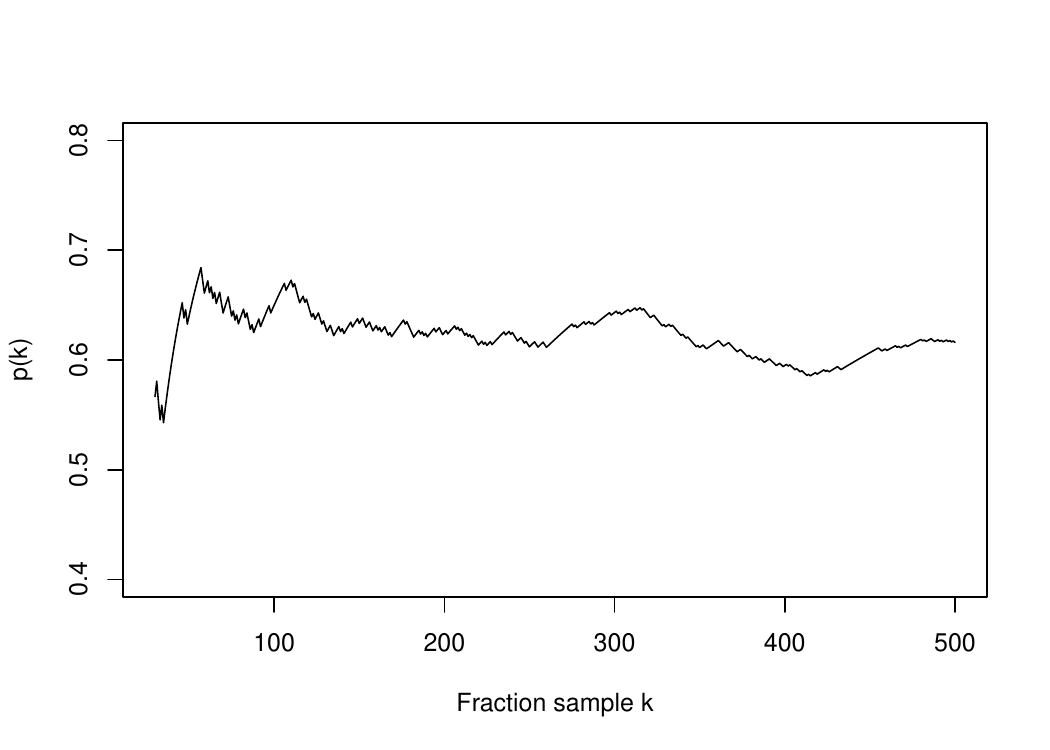}\hfill
\includegraphics[scale=0.42]{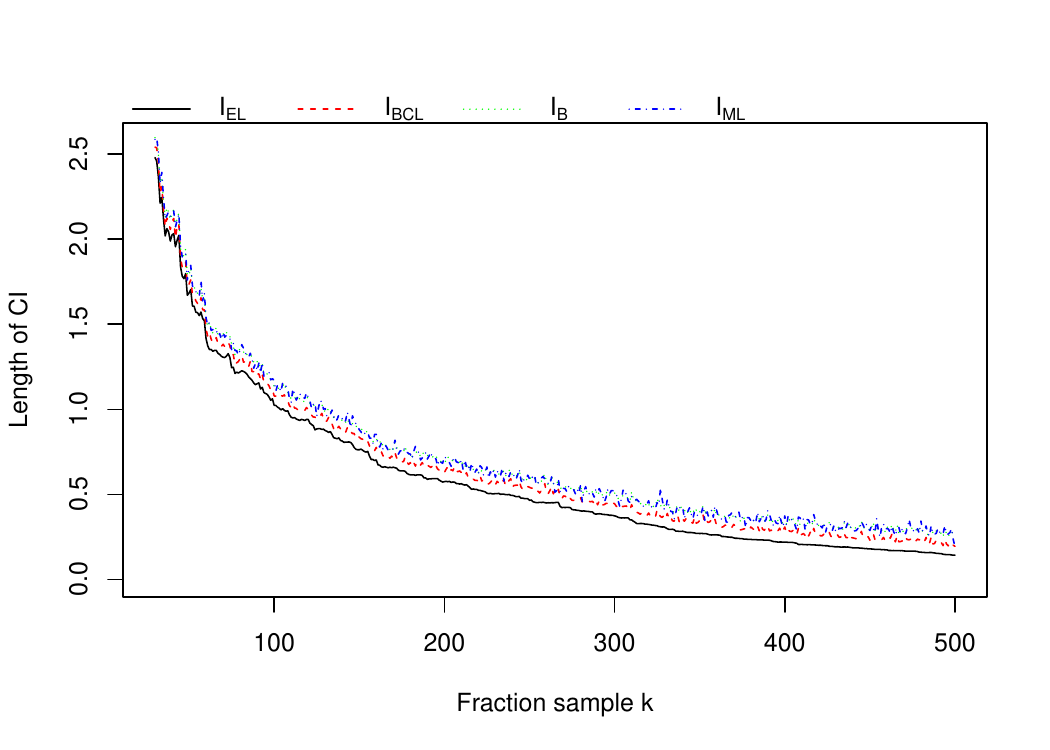}\\
\includegraphics[keepaspectratio=true,scale=0.42]{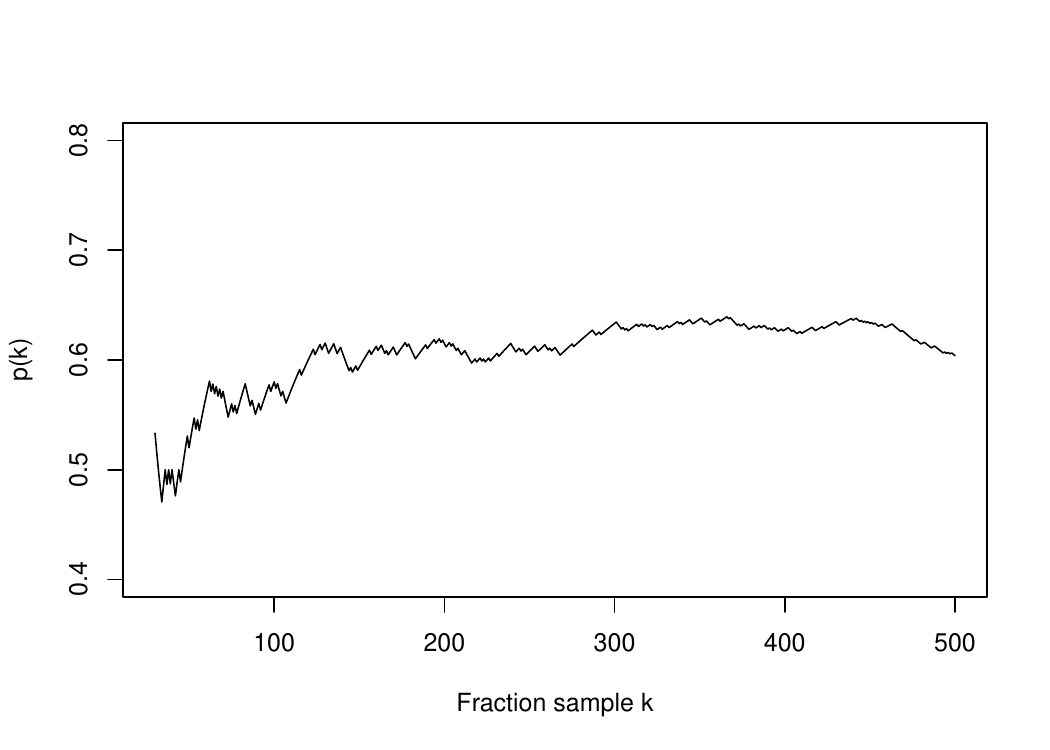}\hfill
\includegraphics[scale=0.42]{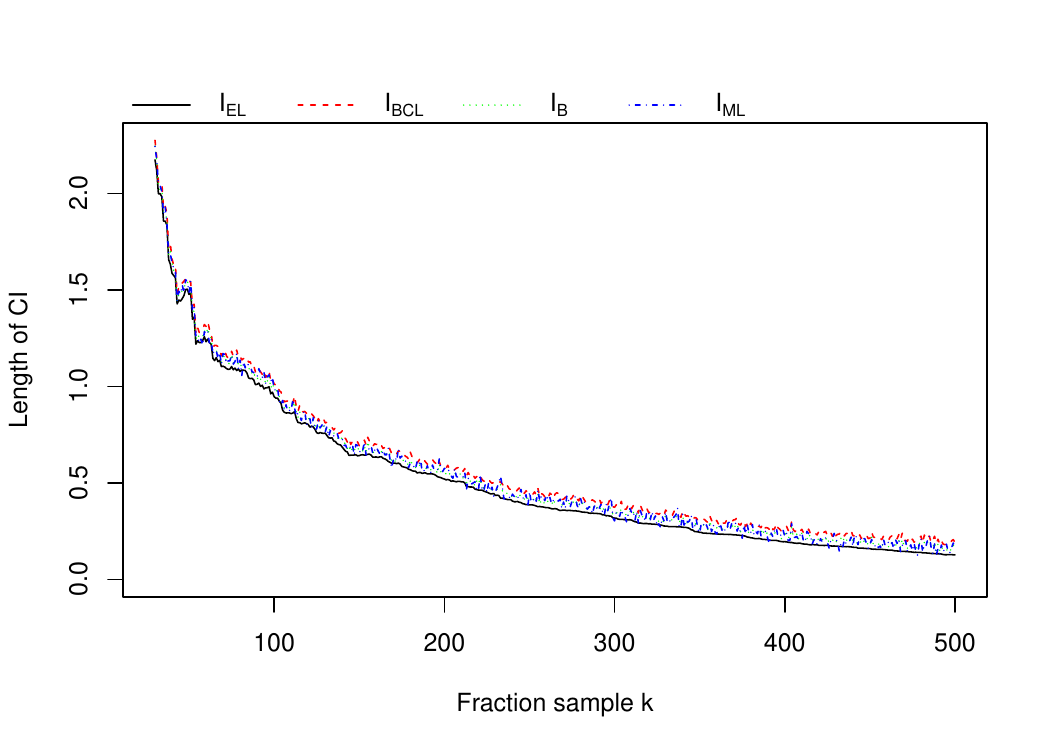}\\

\end{center}
\caption{Coverage probability (left) and average length of confidence region  (right) of $I_{EL}$ (black line), $I_{BCL}$ (red dashed), $I_{B}$ (green dotted),  $I_{ML}$ (blue dotdash) for Australian AIDS survival data ($\mathcal S_1$: first row, $\mathcal S_2$: second row, $\mathcal S_3$: third row).}
\label{figss3}
\end{figure}

Consider the first row of Figure \ref{figss3}. Three distinct phases can be distinguished on the right graph. The first phase, for $k=\{30,\ldots,60\}$ the behaviour of $k\mapsto$ length$(k)$ has a steeper slope, showing a faster decreasing rate. The second phase is characterized with moderate decreasing rate, corresponding to $k$ ranges approximately from 60 to 200. The slope of the graph gradually begins to flatten, indicating a slower deceasing rate (but this phase correspond to the fraction of $k$, where the KS test reject the hypothesis of heavy-tailed distribution). We clearly distinguish the same three phases for all the sub-sets $\mathcal S_1, \mathcal S_2, \mathcal S_3$.  Based on the results in Figure (\ref{figss3}), the length of the confidence region $I_{EL}$ outperforms the other confidence regions and is more preferable.

\section{Conclusion and perspectives}\label{sec6}

In this paper, we investigated several methods for constructing tail index confidence intervals using Bayesian composite likelihood approach, with a particular emphasis on its various adaptations and its utilization within randomly censored data contexts and heavy-tailed distribution is considered. Our investigations suggest that for finite sample sizes, the confidence regions using composite likelihood performs better, both in terms of coverage probabilities and length of the intervals, with a slight superiority for the log-posterior ratio method. This efficiency of  composite likelihood methods can be attributed to its favorable asymptotic properties. Additionally, we touched upon the implementation aspects of Bayesian composite likelihood through a succinct real-data application.

Nonetheless, there remain several pertinent issues deserving consideration. In practical situations, it often arises that some covariate information $W$ is available to the investigator and the distribution of the interest variable $X$ depends on it, which significantly influences and can not be unheeded. In such cases, attention pivots towards estimating the conditional tail index. Adapting our approach to accommodate this context also stands as a subject for our forthcoming research.

\section{Appendix : Proof of theorems and corollary}

\begin{proof}[Sketch of proof of theorem \ref{theo01}]
Let consider $Z_i=\log(X_{n,n-i+1}) - \log(X_{n,n-k})$ and 
\begin{equation*}
 w_{(i)}= \frac{1}{n} \exp \biggl\{ -1-\lambda_1  
 - \lambda_2 \biggl[ \log\left(\frac{X_{n,n-i+1}}{X_{n,n-k}} \right)-\frac{1}{\alpha} \biggr] \biggr\},
\end{equation*}
for $1 \leq i \leq k$, and $q_{(i)}= \frac{1}{n}\exp\{-1-\lambda_1\}$ for $k+1 \leq i \leq n$. Let
$$
A_1(\lambda_1)  = \sum_{i=1}^k q_{(i)} \quad \mbox{and} \quad A_2(\lambda_1)=\sum_{i=1}^k q_{(i)} Z_i
$$
where,
\begin{equation}
    A_1(\lambda_1)=1-\frac{n-k}{n}\exp\{-1-\lambda_1\} \quad \mbox{and} \quad \frac{\sum_{i=1}^k q_{(i)} Z_i}{\sum_{i=1}^k q_{(i)}}=\frac{A_2(\lambda_1)}{A_1(\lambda_1)}.
    \label{proof1.6}
\end{equation}
The equations (\ref{proof1.6}) are equivalent to the constraints (\ref{eq2.10}).

Now, by using the expression of $q_{(i)}$, the second equality in (\ref{proof1.6}) can be written as:
\begin{equation}
    \dfrac{\sum_{i=1}^k \exp \left( -\lambda_2 Z_i\right) Z_i}{\sum_{i=1}^k \exp \left( - \lambda_2 Z_i\right)}= \frac{1}{\alpha}A_1(\lambda_1) -\frac{1}{n}\left[ \frac{1}{\alpha A_1(\lambda_1)} + \frac{\log X_{n,n-k}}{A_1(\lambda_1)}\right].
    \label{proof1.7}
\end{equation}
We establish that a solution ($\lambda_1,\lambda_2$) exists by proving that ($\ast$): for a given $\lambda_1$, it there exists $\lambda_2(\lambda_1)$ such that satisfying (\ref{proof1.7}), and ($\ast\ast$): there exists some $\lambda_1$ such that $\left(\lambda_1,\lambda_2(\lambda_1)\right)$ satisfies $A_1(\lambda_1)=1-\frac{n-k}{n}\exp\{-1-\lambda_1\}$

\textbf{Proof of ($\ast$)}: Let
$$\phi(\lambda)=\dfrac{\sum_{i=1}^k \exp \left( -\lambda_2 Z_i\right) Z_i}{\sum_{i=1}^k \exp \left( - \lambda_2 Z_i\right)}.$$
One can check that $\phi(\lambda) \to Z_1$ when $\lambda \to -\infty$ and $\phi(\lambda) \to Z_k$ when $\lambda \to +\infty$ and $\phi$ is decreasing and continuous on $[Z_1,Z_k]$. Therefore, for a given $\lambda_1$, since $\phi$ is invertible on $[Z_1,Z_k]$, we can define $\lambda_2(\lambda_1)$ as:
\begin{equation}
    \lambda_2(\lambda_1) = \phi^{-1}\left(
\frac{1}{\alpha}A_1(\lambda_1) -\frac{1}{n}\left[ \frac{1}{\alpha A_1(\lambda_1)} + \frac{\log X_{n,n-k}}{A_1(\lambda_1)}\right]
    \right).
    \label{proof1.8}
\end{equation}
It is not difficult to show that $\lambda_2(\lambda_1)$ satisfies (\ref{proof1.7}), and this, $(\lambda_1,\lambda_2(\lambda_1))$ is solution of the second equality in (\ref{proof1.6}).
\textbf{Proof of ($\ast\ast$)}:
To simplify the notation, let $\Tilde{\lambda}=\lambda_2(\lambda_1)$. First, using the Taylor's expansion of $\phi$, 
$$
\phi\left(\pm \frac{1}{\sqrt[4]{k}}\right) = \phi(0) \pm \frac{1}{\sqrt[4]{k}} \phi'(0)\left( 1+o_{\mathbb{P}}(1)\right),
$$
where $\phi(0)=\frac{1}{k}\sum_{i=1}^kZ_i = H(k)$ and 
\begin{equation}
    \phi'(0)=h(k)-\frac{1}{k}\sum_{i=1}^kZ_i^2 = -\frac{1}{\alpha^2}\left(1+ O_{\mathbb{P}}(1) \right),
    \label{proof1.9}
\end{equation}
Using this, we can show that $\mathbb{P}\left(-\frac{1}{\sqrt[4]{k}} < \Tilde{\lambda} < \frac{1}{\sqrt[4]{k}} \right) \to 1$ as $n\to \infty$. A second Taylor's expansion yields
\begin{equation}
    \phi(\Tilde{\lambda})-H(k) = -\Tilde{\lambda} \phi'(0)\left(1+ O_{\mathbb{P}}(k^{-\frac{1}{4}}) \right) = -\Tilde{\lambda} \frac{1}{\alpha^2}\left(1+ O_{\mathbb{P}}(k^{-\frac{1}{4}}) \right),
    \label{proof1.10}
\end{equation}
which implies that $\Tilde{\lambda}=O_{\mathbb{P}}(k^{-\frac{1}{2}})$. Using Taylor's expansion and the results above, we have
\begin{equation}
    \sum_{i=1}^k \exp(-\Tilde{\lambda}Z_i)=k\left( 1+ O_{\mathbb{P}}(k^{-\frac{1}{2}})\right).
    \label{proof1.11}
\end{equation}
Now, let define the function $h(\lambda_1)=A_1(\lambda_1)-\sum_{i=1}^k q_{(i)}$, that we can re-express as function of $\lambda_1$ and $\Tilde{\lambda}$ as:
$$
h(\lambda_1)= A_1(\lambda_1) - \frac{1}{n} \exp\left(-1-\lambda_1 + \frac{1}{\alpha}\Tilde{\lambda}\right) \sum_{i=1}^k \exp(-\Tilde{\lambda}Z_i).
$$
For $\lambda_1^{*}=-1-\log \left( 1+ \frac{\sqrt{k}}{n-k}\right)$ and $\lambda_1^{**}=-1-\log \left( 1- \frac{\sqrt{k}}{n-k}\right)$, it is easy to see that
\begin{equation}
    A_1(\lambda_1^{*}) = \frac{k}{n}\left( 1- \frac{1}{\sqrt{k}}\right) \quad \mbox{ and }  \quad A_1(\lambda_1^{**}) = \frac{k}{n}\left( 1+ \frac{1}{\sqrt{k}}\right).
    \label{proof1.12}
\end{equation}
since:
\begin{itemize}
    \item[\textbf{C1}:] $A_1(\lambda)$ is increasing and continuous function, for every $\lambda_1^{*} \leqslant \lambda \leqslant \lambda_1^{**}$ we have 
    $$
    \frac{k}{n}\left( 1- \frac{1}{\sqrt{k}}\right)\leqslant A_1(\lambda) \leqslant \frac{k}{n}\left( 1+ \frac{1}{\sqrt{k}}\right)
    $$
    \item[\textbf{C2}:] from the asymptotic normality of the Hill estimator we have: $\sqrt{k}\left( H(k)-\alpha^{-1}\right) \xrightarrow{\mathcal{D}} \mathcal{N}(0,\alpha^{-2})$
\end{itemize}
It follows from (\textbf{C1}),  (\textbf{C2}), (\ref{proof1.11}) and (\ref{proof1.12}) that
$$
A_2(\lambda_1^{*}) = -\frac{1}{\alpha}\frac{k}{n}\left( 1- \frac{1}{\sqrt{k}} + o_{\mathbb{P}}(k^{-\frac{1}{2}}) \right), \quad \mbox{ and } \quad A_2(\lambda_1^{**}) = -\frac{1}{\alpha}\frac{k}{n}\left( 1+ \frac{1}{\sqrt{k}} + o_{\mathbb{P}}(k^{-\frac{1}{2}}) \right),
 $$
and
\begin{equation}
    h(\lambda_1^{*})=-\frac{k}{n}\frac{1}{\sqrt{k}}\left(1+o_\mathbb{P}(1)\right) \quad \mbox{ and } \quad h(\lambda_1^{**})=\frac{k}{n}\frac{1}{\sqrt{k}}\left(1+o_\mathbb{P}(1)\right).
    \label{proof1.13}
\end{equation}
for $\lambda_1^{*} \leqslant \lambda \leqslant \lambda_1^{**}$, and by letting $h(\lambda_1^{*})$ and $h(\lambda_1^{**})$ close to 0 as $n\to \infty$, $\lambda_1$ satisfies the first constraint in (\ref{proof1.6}), with $\lambda_2(\lambda_1)$ defined in (\ref{proof1.8}).

Now, for each $\hat{\lambda}_1$ such that $h(\hat{\lambda}_1)=0$, we can show that
\begin{equation}
    1+\hat{\lambda}_1 = O_{\mathbb{P}}\left( \frac{\sqrt{k}}{n}\right).
\end{equation}
Moreover, from (\ref{proof1.10}), we can show that
\begin{equation}
    \Tilde{\lambda}^2 = \left( \dfrac{H(k)-\alpha^-1}{\alpha^{-2}}\right)^2 \left(1+O_{\mathbb{P}}(k^{-1/4}) \right).
    \label{proof1.14}
\end{equation}

Typically as in \cite{peng01}, we have
$$
 \mathbb{L}(\alpha) = (1+o_{\mathbb{P}}(1))\left(n (1+\hat{\lambda}_1)^2 + \Tilde{\lambda}^2 \sum_{i=1}^k  \left( Z_i -H(k)\right)^2\right). 
$$
Finally, by letting $\gamma=\alpha^-1$, it follows from (\ref{proof1.10}), (\ref{proof1.13}) and (\ref{proof1.13}), that
\begin{align*}
    \mathbb{L}(\alpha) &= (1+o_{\mathbb{P}}(1)) O_{\mathbb{P}}(k/n) + (1+o_{\mathbb{P}}(1))\frac{k}{\alpha^2} (1+ O_{\mathbb{P}}(k^{-1/2})) \Tilde{\lambda}^2 \\
    &= o_{\mathbb{P}}(1) + (1+o_{\mathbb{P}}(1)) \left( \sqrt{k} \left( \frac{H(k)- \gamma}{\gamma}\right)^2 \right)\\
    & \xrightarrow{\mathcal{D}} \chi_1^2.
\end{align*}


\end{proof}

\begin{proof}[Sketch of proof of corollary \ref{cor01}]
    To prove the result in corollary (\ref{cor01}), we need to remind some asymptotic results when random censoring is considered:
    By \cite{Ein:01}, $\hat{p}=\frac{1}{k}\sum_{i=1}^k \delta_{[n-k+i]}$  converges in probability to $p:=\frac{\alpha_0}{\alpha_0+\beta_0}$. Thus there exists a subsequence of $k$ along which $\hat{p}$ converges almost surely (a.s.) to $p$. From \cite{DHM:88}, $\frac{1}{k}\sum_{i=1}^{k}\log \left(\frac{Z_{n,n-k+i}}{Z_{n,n-k}}\right)$ converges a.s. to the EVI $\frac{1}{\alpha_0+\beta_0}$ of $Z$ as $n \rightarrow \infty$ and thus, also converges to $\frac{1}{\alpha_0+\beta_0}$ a.s. along the subsequence.
    \cite{amer} showed the consistency of the Bayesian estimators of the tail index under right randomly censored data, and for the same subsequence $k=k_n$, $\sqrt{k}(\alpha - \hat{\alpha}_{MPE}^{(J)}) \xrightarrow{\mathcal{D}} \mathcal{N}(0,\frac{\alpha_0^2}{p})$. Under the corollary assumption and by the Slutsky's theorem, we can write
\begin{align}
\begin{cases}
\frac{1}{k} \sum_{i=1}^k (\frac{k\hat{p}_k}{V_i} - \alpha) &= O_\mathbb{P}\left(k^{-1/2}+|A(n/k)|\right)\\
\frac{1}{k} \sum_{i=1}^k \left(\frac{k\hat{p}_k}{V_i} - \alpha\right)^2 &= \frac{p}{\alpha^2} + O_\mathbb{P}\left(k^{-1/2}+|A(n/k)|\right). \label{eq.proof3.1}
\end{cases}
\end{align}
        
    Since the weights $w_i$ appearing in the $\ell(\alpha)$ are given by 
\begin{equation*}
w_i = k^{-1}\{1 + \lambda(\frac{k\hat{p}_k}{V_i} - \alpha)\}^{-1},
\end{equation*}
we have $1 + \lambda(\frac{k\hat{p}_k}{V_i} - \alpha) \geq 0$, and thus
\begin{equation*}
|1 + \lambda(\frac{k\hat{p}_k}{V_i} - \alpha)|^{-1} \geq \{1 + |\lambda| \max_i |\frac{k\hat{p}_k}{V_i} - \alpha|\}^{-1}.
\end{equation*}
By definition of $\lambda$ we have
\begin{align*}
0 &= k^{-1}\left|\sum_{i=1}^k \frac{\lambda(\frac{k\hat{p}_k}{V_i} - \alpha)^2}{1 + \lambda(\frac{k\hat{p}_k}{V_i} - \alpha)} - (\frac{k\hat{p}_k}{V_i} - \alpha)\right|\\
 &\geq k^{-1} \sum_{i=1}^k \frac{|\lambda|(\frac{k\hat{p}_k}{V_i} - \alpha)^2}{|1 + \lambda(\frac{k\hat{p}_k}{V_i} - \alpha)|} - \left| \frac{1}{k} \sum_{i=1}^k (\frac{k\hat{p}_k}{V_i} - \alpha)\right| \\
 &\geq  k^{-1} \sum_{i=1}^k \frac{|\lambda|(\frac{k\hat{p}_k}{V_i} - \alpha)^2}{\{1 + |\lambda| \max_{1 \leq i \leq k} |\frac{k\hat{p}_k}{V_i} - \alpha|\}} - \left| \frac{1}{k} \sum_{i=1}^k (\frac{k\hat{p}_k}{V_i} - \alpha)\right|,\\
\end{align*}
which implies that
\begin{equation*}
|\lambda| \left[ \frac{1}{k} \sum_{i=1}^k (\frac{k\hat{p}_k}{V_i} - \alpha)^2 -  \max_{1 \leq i \leq k} \left|\frac{k\hat{p}_k}{V_i} - \alpha\right| \left| \frac{1}{k} \sum_{i=1}^k (\frac{k\hat{p}_k}{V_i} - \alpha) \right|\right] \leq \left| \frac{1}{k} \sum_{i=1}^k (\frac{k\hat{p}_k}{V_i} - \alpha) \right|
\end{equation*}
Hence, by (\ref{eq.proof3.1}), $\lambda = O_\mathbb{P}\left(k^{-1/2}+|A(n/k)|\right)$. Furthermore, by Taylor's expanding of (\ref{eq.lambda}) we have
\begin{align*}
0 &= \frac{1}{k} \sum_{i=1}^k (\frac{k\hat{p}_k}{V_i} - \alpha)\{1 - \lambda(\frac{k\hat{p}_k}{V_i} - \alpha) - \lambda^2(\frac{k\hat{p}_k}{V_i} - \alpha)^2 - \lambda^3(\frac{k\hat{p}_k}{V_i} - \alpha)^3 - \cdots\}\\
 &= \frac{1}{k} \sum_{i=1}^k (\frac{k\hat{p}_k}{V_i} - \alpha) - \lambda \frac{1}{k} \sum_{i=1}^k (\frac{k\hat{p}_k}{V_i} - \alpha)^2  + O_\mathbb{P}\left((k^{-1/2}+|A(n/k)|)^2\right),
\end{align*}
i.e.
\begin{equation}
\lambda = \frac{k^{-1} \sum_{i=1}^k (\frac{k\hat{p}_k}{V_i} - \alpha)}{k^{-1} \sum_{i=1}^k (\frac{k\hat{p}_k}{V_i} - \alpha)^2} + O_\mathbb{P}\left((k^{-1/2}+|A(n/k)|)^2\right). \label{eq.proof.3}
\end{equation}
Thus, by(\ref{eq.proof3.1}) and (\ref{eq.proof.3}) and the assumption $\sqrt{k}A(n/k) \to 0$, we have, as $n \to \infty$,
\begin{align*}
\ell(\alpha) &= 2 \sum_{i=1}^k \log\{1 + \lambda(\frac{k\hat{p}_k}{V_i} - \alpha)\}\\
       &= 2\lambda \sum_{i=1}^k (\frac{k\hat{p}_k}{V_i} - \alpha) - \lambda^2 \sum_{i=1}^k (\frac{k\hat{p}_k}{V_i} - \alpha)^2 + O_\mathbb{P}\left((k^{-1/2}+|A(n/k)|)^3\right)\\
       &= k \dfrac{\left[k^{-1} \sum_{i=1}^k (\frac{k\hat{p}_k}{V_i} - \alpha)\right]^2}{k^{-1} \sum_{i=1}^k (\frac{k\hat{p}_k}{V_i} - \alpha)^2} + O_\mathbb{P}\left(k(k^{-1/2}+|A(n/k)|)^3\right)\\
       &= k \dfrac{\left[k^{-1} \sum_{i=1}^k (\frac{k\hat{p}_k}{V_i} - \alpha) + O_\mathbb{P}\left(|A(n/k)|\right)\right]^2 }{k^{-1} \sum_{i=1}^k (\frac{k\hat{p}_k}{V_i} - \alpha)^2 + O_\mathbb{P}\left(|A(n/k)|\right)} + o_\mathbb{P}\left(1\right)\\
       &= k \dfrac{\left[k^{-1} \sum_{i=1}^k (\frac{k\hat{p}_k}{V_i} - \alpha)  \right]^2 }{k^{-1} \sum_{i=1}^k (\frac{k\hat{p}_k}{V_i} - \alpha)^2 } + O_\mathbb{P}\left(\sqrt{k}|A(n/k)|\right) + o_\mathbb{P}\left(1\right)\\
       &\overset{d}{\to} \chi^2_1.
\end{align*}
\end{proof}

\bibliographystyle{apalike}
\bibliography{references}

\begin{thebibliography}{}

\bibitem[Ameraoui et~al., 2016]{amer}
Ameraoui, A., Boukhetala, K., and Dupuy, J. (2016).
\newblock Bayesian estimation of the tail index of a heavy tailed distribution under random censoring.
\newblock {\em Computational Statistics \& Data Analysis}, 104:148--168.
\newblock \url{https://doi.org/10.1016/j.csda.2016.06.009}.

\bibitem[Baysal and Staum, 2008]{bay01}
Baysal, R. and Staum, J. (2008).
\newblock Empirical likelihood for value-at-risk and expected shortfall.
\newblock {\em Journal of Risk}, 11(1):3--32.
\newblock \url{https://doi.org/10.21314/JOR.2008.185}.

\bibitem[Beirlant et~al., 2004]{Bei:04}
Beirlant, J., Goegebeur, Y., Teugels, J., and Segers, J. (2004).
\newblock {\em Statistics of Extremes}.
\newblock John Wiley \& Sons Ltd, Chichester.
\newblock ISBN: 978-0-471-97647-9. \url{https://doi.org/10.1002/0470012382}.

\bibitem[Beirlant et~al., 2007]{Bei:07}
Beirlant, J., Guillou, A., Dierckx, G., and Fils-Villetard, A. (2007).
\newblock Estimation of the extreme value index and extreme quantiles under random censoring.
\newblock {\em Extremes}, 10(3):151--174.
\newblock \url{https://doi.org/10.1007/s10687-007-0039-x}.

\bibitem[Brahimi et~al., 2015]{Bra:13}
Brahimi, B., Meraghni, D., and Necir, A. (2015).
\newblock Gaussian approximation to the extreme value index estimator of a heavy-tailed distribution under random censoring.
\newblock {\em Mathematical Methods of Statistics}, 24(4):266--279.
\newblock \url{https://doi.org/10.3103/S106653071504002X}.

\bibitem[Cabras and Castellanos, 2011]{Cab:11}
Cabras, S. and Castellanos, M. (2011).
\newblock A bayesian approach for estimating extreme quantiles under a semiparametric mixture model.
\newblock {\em Astin Bulletin}, 41(1):87--106.
\newblock \url{https://doi.org/10.2143/AST.41.1.2084387}.

\bibitem[Casella and Berger, 2024]{Cas:24}
Casella, G. and Berger, R. (2024).
\newblock {\em Statistical Inference (2nd ed.).}
\newblock Chapman and Hall/CRC.
\newblock ISBN: 978-1-003-45628-5. \url{https://doi.org/10.1201/9781003456285}.

\bibitem[Chan and So, 2017]{chan:17}
Chan, R. and So, M. (2017).
\newblock On the performance of the bayesian composite likelihood estimation of max-stable processes.
\newblock {\em Journal of Statistical Computation and Simulation}, 87(15):2869--2881.
\newblock \url{https://doi.org/10.1080/00949655.2017.1342824}.

\bibitem[Coles and Powell, 1996]{Coles:96}
Coles, S. and Powell, E. (1996).
\newblock Bayesian methods in extreme value modelling: a review and new developments.
\newblock {\em International Statistical Review}, 64(1):119--136.
\newblock \url{https://doi.org/10.2307/1403426}.

\bibitem[de~Haan and Ferreira, 2006]{dHan01}
de~Haan, L. and Ferreira, A. (2006).
\newblock {\em Extreme Value Theory}.
\newblock Springer, New York.
\newblock ISBN: 978-0-387-23946-0. \url{https://doi.org/10.1007/0-387-34471-3}.

\bibitem[de~Zea~Bermudez and Kotz, 2010a]{Ber:02}
de~Zea~Bermudez, P. and Kotz, S. (2010a).
\newblock Parameter estimation of the generalized pareto distribution-part i.
\newblock {\em Journal of Statistical Planning and Inference,}, 140(6):1353--1373.
\newblock \url{https://doi.org/10.1016/j.jspi.2008.11.019}.

\bibitem[de~Zea~Bermudez and Kotz, 2010b]{Ber:01}
de~Zea~Bermudez, P. and Kotz, S. (2010b).
\newblock Parameter estimation of the generalized pareto distribution-part ii.
\newblock {\em Journal of Statistical Planning and Inference,}, 140(6):1374--1388.
\newblock \url{https://doi.org/10.1016/j.jspi.2008.11.020}.

\bibitem[Deheuvels et~al., 1988]{DHM:88}
Deheuvels, P., Haeusler, E., and Mason, M. (1988).
\newblock Almost sure convergence of the hill estimator.
\newblock {\em Mathematical Proceedings of the Cambridge Philosophical Society}, 104(2):371--381.
\newblock \url{https://doi.org/10.1017/S0305004100065531}.

\bibitem[Diebolt et~al., 2005]{Dieb:05}
Diebolt, J., El-Aroui, M., Garrido, M., and Girard, S. (2005).
\newblock Quasi-conjugate bayes estimates for gpd parameters and application to heavy tails modelling.
\newblock {\em Extremes}, 8(1):57--78.
\newblock \url{https://doi.org/10.1007/s10687-005-4860-9}.

\bibitem[do~Nascimento et~al., 2012]{Nasc:12}
do~Nascimento, F., Gamerman, D., and Lopes, H. (2012).
\newblock A semiparametric bayesian approach to extreme value estimation.
\newblock {\em Statistics and Computing}, 22:661--675.
\newblock \url{https://doi.org/10.1007/s11222-011-9270-z}.

\bibitem[Efron, 1993]{Efr:93}
Efron, B. (1993).
\newblock Bayes and likelihood calculations from confidence intervals.
\newblock {\em Biometrika}, 80:3--26.
\newblock \url{https://doi.org/10.1093/biomet/80.1.3}.

\bibitem[Einmahl et~al., 2008]{Ein:01}
Einmahl, J., Fils-Villetard, A., and Guillou, A. (2008).
\newblock Statistics of extremes under random censoring.
\newblock {\em Bernoulli}, 14(1):207--227.
\newblock \url{https://doi.org/10.3150/07-BEJ104}.

\bibitem[Embrechts et~al., 1997]{Emb:01}
Embrechts, P., Kluppelberg, C., and Mikosch, T. (1997).
\newblock {\em Modelling Extremal Events For Insurance and Finance}.
\newblock Springer-Verlag, Berlin.
\newblock ISBN: 978-3-540-609315. \url{https://doi.org/10.1007/978-3-642-33483-2}.

\bibitem[Gomes and Neves, 2011]{GN:11}
Gomes, M. and Neves, M. (2011).
\newblock Estimation of the extreme value index for randomly censored data.
\newblock {\em Biometrical Letters}, 48(1):1--22.

\bibitem[He et~al., 2016]{he01}
He, S., Liang, W., Shen, J., and Yang, G. (2016).
\newblock Empirical likelihood for right censored lifetime data.
\newblock {\em Journal of the American Statistical Association}, 111(514):646--655.
\newblock \url{https://doi.org/10.1080/01621459.2015.1024058}.

\bibitem[Hill, 1975]{hill01}
Hill, B. (1975).
\newblock A simple general approach to inference about the tail of a distribution.
\newblock {\em The Annals of Statistics}, 3(5):1163--1174.
\newblock \url{https://doi.org/10.1214/aos/1176343247}.

\bibitem[Jeffreys, 1998]{Jeff:61}
Jeffreys, H. (1998).
\newblock {\em Theory of Probability.}
\newblock Oxford University Press.
\newblock ISBN: 978-0-198503682. \url{https://doi.org/10.1093/oso/9780198503682.001.0001}.

\bibitem[Koning and Peng, 2008]{Kon:08}
Koning, A. and Peng, L. (2008).
\newblock Goodness-of-fit tests for a heavy tailed distribution.
\newblock {\em Journal of Statistical Planning and Inference}, 138(12):3960--3981.
\newblock \url{https://doi.org/10.1016/j.jspi.2008.02.013}.

\bibitem[Lazar, 2003]{laz01}
Lazar, N. (2003).
\newblock Bayesian empirical likelihood.
\newblock {\em Biometrika}, 90(2):319--326.
\newblock \url{https://www.jstor.org/stable/30042042}.

\bibitem[Li and Qi, 2019]{li01}
Li, Y. and Qi, Y. (2019).
\newblock Adjusted empirical likelihood method for the tail index of a heavy-tailed distribution.
\newblock {\em Statistics \& Probability Letters}, 152:50--58.
\newblock \url{https://doi.org/10.1016/j.spl.2019.04.015}.

\bibitem[Liang and Dai, 2021]{liang01}
Liang, W. and Dai, H. (2021).
\newblock Empirical likelihood based on synthetic right censored data.
\newblock {\em Statistics \& Probability Letters}, 169:0167--7152.
\newblock \url{https://doi.org/10.1016/j.spl.2020.108962}.

\bibitem[Lindsay, 1988]{Lind:88}
Lindsay, B. (1988).
\newblock Composite likelihood methods.
\newblock {\em Contemporary Mathematics}, 80:221--239.
\newblock \url{http://dx.doi.org/10.1090/conm/080/999014}.

\bibitem[Lu and Peng, 2002]{lu02}
Lu, J. and Peng, L. (2002).
\newblock Likelihood based confidence intervals for the tail index.
\newblock {\em Extremes}, 5:337--352.
\newblock \url{https://doi.org/10.1023/A:1025163807024}.

\bibitem[Owen, 1988]{owen01}
Owen, A.~B. (1988).
\newblock Empirical likelihood ratio confidence intervals for a single functional.
\newblock {\em Biometrika}, 75(2):237--249.
\newblock \url{https://doi.org/10.2307/2336172}.

\bibitem[Owen, 1990]{owen02}
Owen, A.~B. (1990).
\newblock Empirical likelihood ratio confidence regions.
\newblock {\em The Annals of Statistics}, 18(1):90--–120.
\newblock \url{https://doi.org/10.1214/aos/1176347494}.

\bibitem[Owen, 2001]{owen03}
Owen, A.~B. (2001).
\newblock {\em Empirical Likelihood}.
\newblock Chapman and Hall/CRC, London.
\newblock ISBN: 978-0-429128998. \url{https://doi.org/10.1201/9781420036152}.

\bibitem[Owen, 2013]{owen04}
Owen, A.~B. (2013).
\newblock Self-concordance for empirical likelihood.
\newblock {\em The Canadian Journal of Statistics}, 41(3):387--397.
\newblock \url{https://doi.org/10.1002/cjs.11183}.

\bibitem[Pauli et~al., 2011]{Pau:12}
Pauli, F., Racugno, W., and Ventura, L. (2011).
\newblock Bayesian composite marginal likelihoods.
\newblock {\em Statistica Sininca}, 21:149--164.
\newblock \url{http://www.jstor.org/stable/24309266}.

\bibitem[Peng and Qi, 2006]{peng01}
Peng, L. and Qi, Y. (2006).
\newblock Confidence region for high quantiles of a heavy tailed distribution.
\newblock {\em The Annals of Statistics}, 34(4):1964--1986.
\newblock \url{https://doi.org/10.1214/009053606000000416}.

\bibitem[Peng et~al., 2015]{peng02}
Peng, L., Wang, X., and Zheng, Y. (2015).
\newblock Empirical likelihood inference for haezendonck-goovaerts risk measure.
\newblock {\em European Actuarial Journal}, 5(2):427--445.
\newblock \url{https://doi.org/10.1007/s13385-015-0113-8}.

\bibitem[Qi, 2008]{qi01}
Qi, Y. (2008).
\newblock Bootstrap and empirical likelihood methods in extremes.
\newblock {\em Extremes}, 11:81--97.
\newblock \url{https://doi.org/10.1007/s10687-700-0049-8}.

\bibitem[Ribatet et~al., 2012]{Rib:12}
Ribatet, M., Cooley, D., and Davison, A. (2012).
\newblock Bayesian inference from composite likelihoods, with an application to spatial extremes.
\newblock {\em Statistica Sinica}, 22(2012):813--845.
\newblock \url{http://dx.doi.org/10.5705/ss.2009.248}.

\bibitem[Ritchie and Rosado, 2022]{rit:22}
Ritchie, H. and Rosado, P. (2022).
\newblock Natural disasters.
\newblock {\em Our World in Data}.
\newblock \url{https://ourworldindata.org/natural-disasters}.

\bibitem[Stein, 2023]{Stn:23}
Stein, M. (2023).
\newblock A weighted composite log-likelihood approach to parametric estimation of the extreme quantiles of a distribution.
\newblock {\em Extremes}, 26(3):469--507.
\newblock \url{https://doi.org/10.1007/s10687-023-00466-w}.

\bibitem[Team, 2008]{softR}
Team, R. D.~C. (2008).
\newblock {\em R: A language and environment for statistical computing.}
\newblock R Foundation for Statistical Computing, Vienna, Austria.
\newblock 3-900051-07-0. \url{http://www.R-project.org}.

\bibitem[Venables and Ripley, 2002]{VenRip:02}
Venables, W. and Ripley, B. (2002).
\newblock {\em Modern Applied Statistics with S.}
\newblock Springer, New York.
\newblock ISBN: 978-0-387-95457-8. \url{https://doi.org/10.1007/978-0-387-21706-2}.

\bibitem[Weissman, 1978]{weiss01}
Weissman, I. (1978).
\newblock Estimation of parameters and large quantiles based on the k largest observations.
\newblock {\em Journal of the American Statistical Association}, 73(364):812--815.
\newblock \url{https://doi.org/10.1080/01621459.1978.10480104}.

\bibitem[Worms and Worms, 2014]{Wor:14}
Worms, J. and Worms, R. (2014).
\newblock New estimators of the extreme value index under random right censoring, for heavy tailed distributions.
\newblock {\em Extremes}, 17(2):337--358.
\newblock \url{https://doi.org/10.1007/s10687-014-0189-6}.

\bibitem[Xi and Keilegom, 2006]{chen01}
Xi, C. and Keilegom, I.~V. (2006).
\newblock A review on empirical likelihood methods for regression.
\newblock {\em Test}, 18(3):415--447.
\newblock \url{https://doi.org/10.1007/s11749-009-0159-5}.

\bibitem[Yan and Zhang, 2017]{zhen01}
Yan, Z. and Zhang, J. (2017).
\newblock Adjusted empirical likelihood for value-at-risk and expected shortfall.
\newblock {\em Communications in Statistics - Theory and Methods}, 46(5):2580--2591.
\newblock \url{https://doi.org/10.1080/03610926.2014.1002933}.

\end{thebibliography}

\end{document}